\documentclass[10pt]{article}
\usepackage{amsmath, amssymb, amsfonts, amsthm, epsfig, float, graphicx, sectsty}
\usepackage[all]{xy}
\title{Uniform Diameter Bounds in Branch Groups}
\author{Henry Bradford}

\newcommand{\Addresses}{{
  \bigskip
  \footnotesize

H. Bradford, \textsc{Universit\"{a}t G\"{o}ttingen, Germany}\par\nopagebreak
  \textit{E-mail address:} \texttt{henry.bradford@mathematik.uni-goettingen.de}

}}

\newtheorem{thm}{Theorem}[section]
\newtheorem{lem}[thm]{Lemma}
\newtheorem{propn}[thm]{Proposition}
\newtheorem{coroll}[thm]{Corollary}
\newtheorem{defn}[thm]{Definition}
\newtheorem{ex}[thm]{Example}

\newtheorem{rmrk}[thm]{Remark}

\DeclareMathOperator{\Alt}{Alt}
\DeclareMathOperator{\Aut}{Aut}
\DeclareMathOperator{\Cay}{Cay}
\DeclareMathOperator{\diam}{diam}

\DeclareMathOperator{\SL}{SL}
\DeclareMathOperator{\Syl}{Syl}
\DeclareMathOperator{\Stab}{Stab}
\DeclareMathOperator{\Sym}{Sym}

\begin{document}

\maketitle

\begin{abstract}
Let $G$ be either the Grigorchuk $2$-group 
or one of the Gupta-Sidki $p$-groups. 
We give new upper bounds for the diameters of the quotients of $G$ by its level stabilisers, 
as well as other natural sequences of finite-index normal subgroups. 
Our bounds are independent of the generating set, 
and are polylogarithmic functions of the group order, with explicit degree. 
Our proofs utilize a version of the profinite 
Solovay-Kitaev procedure, the branch structure of $G$, 
and in certain cases, existing computations of the lower central series of $G$. 
\end{abstract}

\section{Introduction}

Let $G$ be a finite group, and $S \subseteq G$ be a generating set. 
The \emph{diameter of $G$ with respect to $S$} is defined to be:
\begin{center}
$\diam (G,S) = \min \lbrace n \in \mathbb{N} : B_S (n) = G \rbrace$, 
\end{center}
where $B_S (n)$ is the (closed) ball of radius $n$ about the identity in the word metric defined by $S$ on $G$. The \emph{diameter of $G$}, denoted $\diam(G)$, is then the maximal value of $\diam (G,S)$ as $S$ ranges over all generating subsets of $G$. 
In this paper we give upper bounds for the diameters of natural families of finite quotients of certain \emph{branch groups}. 

\subsection{Statement of Results}

\begin{thm} \label{stabdiamcorollGrig}
Let $\mathfrak{G}$ be the Grigorchuk $2$-group. Then: 
\begin{align*}
\diam (\mathfrak{G} / \Stab_{\mathfrak{G}} (n)) 
& = O \big( \exp(\log (35)n) \big) \\
& = O\big( \log \lvert \mathfrak{G} : \Stab_{\mathfrak{G}} (n) \rvert^{\log (35) /\log (2)}\big)\text{.}
\end{align*}
\end{thm}

We shall define the sequence of \emph{level stabilisers} 
$\Stab_G (n)$ for a group $G$ acting on a rooted tree in 
subsection \ref{branchsection}. 
Our proof makes extensive use of the description of the lower central series $(\gamma_n (\mathfrak{G}))_n$ of $\mathfrak{G}$ given in \cite{BarGri}, building on \cite{Roz} (see also \cite{Bar}). The results of these papers  facilitate an explicit description of the restriction to $\gamma_n (\mathfrak{G})$ 
of the action of $\mathfrak{G}$ 
on the binary rooted tree. Indeed, 
Theorem \ref{stabdiamcorollGrig} 
is proved as a consequence of the following. 

\begin{thm} \label{mainthmGrig}
\begin{align*}
\diam (\mathfrak{G} / \gamma_n (\mathfrak{G})) 
& =O \big( n^{\log (35) / \log (2)} \big) \\
& =O \big( \log \lvert \mathfrak{G} : \gamma_n (\mathfrak{G}) \rvert^{\log (35) / \log (2)} \big)\text{.}
\end{align*}
\end{thm}

Recall that $\log (35) / \log (2) \approx 5.129$. 
We now turn to our results on the Gupta-Sidki $p$-groups. 

\begin{thm} \label{stabcorollGS}
Let $p$ be an odd prime. 
Let $\Gamma_{(p)}$ be the Gupta-Sidki $p$-group. Then: 
\begin{align*}
\diam (\Gamma_{(p)} / \Stab_{\Gamma_{(p)}} (n)) 
& = O_p \big( \exp (\log (C_p)n) \big) \\
& = O_p \big( \log \lvert \Gamma_{(p)} : \Stab_{\Gamma_{(p)}} (n) \rvert^{\log (C_p)/\log(p)} \big)
\end{align*}
where $C_p = 3 \cdot 4^p - 2^p (p+8) + 7$. 
\end{thm}

Theorem \ref{stabcorollGS} is a consequence of our next result. 
Let $K$ be the derived subgroup of $\Gamma_{(p)}$. 

\begin{thm} \label{mainthmGS}
Let $C_p$ is as in Theorem \ref{stabcorollGS}. 
\begin{align*}
\diam (\Gamma_{(p)} / K^{(\times p^n)}) 
& = O_p \big( \exp (\log (C_p)n) \big) \\
& = O_p \big(\log \lvert \Gamma_{(p)} : K^{(\times p^n)} \rvert^{\log (C_p)/\log(p)}\big)\text{.}
\end{align*}

\end{thm}

Here $K^{(\times p^n)}$ denotes the Cartesian product of $p^n$ 
copies of $K$. 
We define the natural embeddings of the $K^{(\times p^n)}$ 
as finite-index normal subgroups of $\Gamma_{(p)}$ 
in subsection \ref{branchsection}. 
For now suffice to say that there are inclusions 
$K^{(\times p^n)} \leq \Stab_{\Gamma_{(p)}} (n+1)$
and that given certain well-known bounds on the orders of the relevant groups, 
Theorem \ref{stabcorollGS} quickly follows from these inclusions and 
Theorem \ref{mainthmGS}. 
In the case $p=3$, 
we can exploit the description of the lower central series of 
$\Gamma_{(3)}$ 
given by Bartholdi \cite{Bar} to also deduce the following. 

\begin{thm} \label{centralcorollGS}
\begin{align*}
\diam (\Gamma_{(3)} / \gamma_n (\Gamma_{(3)})) 
& = O \big(  n^{\log (111)/\log(1+\sqrt{2})}  \big) \\
& = O \big(\log \lvert \Gamma_{(3)} : \gamma_n (\Gamma_{(3)}) \rvert^{\log (111) / \log (3)} \big) \text{.}
\end{align*}
\end{thm}

For the Gupta-Sidki $3$-group therefore, our bounds for the diameter are  polylogarithmic in the order of the group, with degree $\log (C_3) / \log (3)$ \linebreak$= \log (111) / \log (3) \approx 4.287$. As $p$ grows, 
the degree of the polylogarithm grows, proportional to $p/\log(p)$. 
The implied constants in Theorems \ref{stabdiamcorollGrig}-\ref{centralcorollGS} 
may all be explicitly computed from our proofs. 

$\mathfrak{G} / \Stab_{\mathfrak{G}} (n)$ and 
$\Gamma_{(p)} / \Stab_{\Gamma_{(p)}} (n)$ are transitive 
imprimitive permutation groups on, respectively, 
$2^n$ and $p^n$ points. 
As such, Theorems \ref{stabdiamcorollGrig} and \ref{stabcorollGS} 
provide new examples of transitive subgroups of $\Sym(N)$ 
($N$ a power of a fixed prime) whose diameters are polynomially bounded in $N$.  
In the next subsection we will further contextualise Theorems \ref{stabdiamcorollGrig} and \ref{stabcorollGS} 
within the existing literature on diameters of permutation groups. 
$\mathfrak{G}$ and $\Gamma_{(p)}$ are particularly famous members of 
the class of branch groups, and have been extensively studied since 
their introduction, respectively in \cite{Grig1} and \cite{GupSid}. 
They will be defined precisely in subsection \ref{branchsection}. 

\subsection{Background and Structure of the Paper} \label{backgroundsubsect}

$\mathfrak{G}$ is one of the most exotic objects in geometric group theory: it is a finitely generated infinite $2$-group, so provides a counterexample to the General Burnside Problem; it is a group of \emph{intermediate growth}, and indeed was the first example of such a group to be constructed; it is amenable but not elementary amenable; 
it is a residually finite just-infinite group; 
it admits no faithful representation over any field, and 
every finite $2$-group embeds into it as a subgroup \cite{dlHa}. 
Similarly, $\Gamma_{(p)}$ is a finitely generated infinite $p$-group; 
contains a copy of every finite $p$-group, 
and shares many of the other aforementioned properties 
with $\mathfrak{G}$, 
though its growth is the subject of ongoing discussion. 

An understanding of the broader class of \emph{branch groups}, 
to which $\mathfrak{G}$ and $\Gamma_{(p)}$ belong, 
has now become a crucial part of the toolkit of the modern 
geometric or profinite group theorist. 
In one sense this is no surprise, bearing in mind Wilson's 
classification of just-infinite groups \cite{Wilson}, 
in which branch groups comprise an important case 
(although historically, this may only be noted in hindsight, 
since the class of branch groups was not formally defined until 
some time after Wilson's theorem \cite{Grig3}). 
What was perhaps less expected was the extraordinarily 
relevance branch groups would prove to have, 
to subjects as diverse as (but by no means limited to) 
decision problems, finite automata, 
spectral graph theory, fractal spaces, 
and exotic phenomena in the domains of word growth, 
subgroup growth and other asymptotic invariants of infinite groups
\cite{BaGrSu}. 

Seperate from, but roughly concurrent with these developments, 
interest was growing in generation problems in various 
families of groups, including the search for good diameter bounds. 
From the start, particular attention was paid to upper bounds for permutation groups, owing to connections with problems 
in theoretical computer science.  
These included membership testing protocols in 
computational group theory and complexity analysis of deciding solvability of various combinatorial puzzles \cite{DriFur, KoMiSp}, 
the most famous of which is of course the Rubick's cube. 
More recently intense interest in diameters of finite groups 
has been renewed, motivated by connections with 
such diverse topics as expander graphs, 
approximate groups, the Banach-Ruziewicz problem, 
Apollonian circle packings, 
sum-product phenomena in fields and affine sieves. 
The modern study of diameters of groups 
is therefore an extremely rich and diverse subject, 
and one which we cannot hope to fully capture here; 
we instead refer the interested reader to \cite{Helf,Tao} 
and references therein for an overview of the recent developments. 

Among the key tools to have featured in proofs of good upper bounds for 
the diameters of finite groups is the \emph{Solovay-Kitaev procedure}. 
This is a method which was originally used in the context of 
compact complex Lie groups (where it had applications 
to problems in quantum computer science) 
but which translated readily to the setting of abstract 
or profinite groups. Given a group $\Gamma$ and a descending 
sequence $(\Gamma_i)_i$ of finite-index normal subgroups, 
the Solovay-Kitaev procedure proves, 
under suitable additional assumptions, 
an upper bound on the diameters 
of the finite quotient groups $\Gamma / \Gamma_i$ 
by induction on the sequence $(\Gamma_i)_i$. 
A crucial ingredient facilitating the induction step 
is that any element of a later term in the sequence 
should be expressible as 
(or at least sufficiently closely approximable by) 
the product of a small number of commutators 
of elements lying in earlier terms. 
Exactly what this means in practice 
will depend on the specific group $\Gamma$ and sequence of 
subgroups with which we are working. 

In the above setting the Solovay-Kitaev procedure 
was first used by Gamburd and Shahshahani \cite{GambShah} 
to give a polylogarithmic upper bound on the diameter of 
$\SL_2 (\mathbb{Z}/p^n\mathbb{Z})$. 
Their work was subsequently extended, 
first by Dinai \cite{Dinai} to groups of 
$(\mathbb{Z}/p^n\mathbb{Z})$-points of other Chevalley groups, 
then by the author \cite{Brad} to congruence quotients 
of other $p$-adic analytic groups, 
$\mathbb{F}_q [[t]]$-analytic groups 
and the Nottingham groups of finite fields. 

Producing the commutator expressions required 
to facilitate the induction step 
in the Solovay-Kitaev procedure 
requires fairly explicit 
computations of commutators in the subgroups $\Gamma_i$. 
As such, it is very useful in practice 
for the sequence $(\Gamma_i)_i$ 
to be highly ``recurrent'' in some sense, 
so that calculations carried out 
at one level may be translated to others. 
In many of the examples considered in \cite{Brad}, for instance, 
the group $\Gamma$ has a naturally associated Lie algebra 
over a non-archimedean field, such that the $\Gamma_i$ 
may be identified with a descending sequence of balls in the Lie algebra. Under this identification, computations of commutators 
at different levels are simply ``rescalings'' of each other. 

The regular branch structure of $\mathfrak{G}$ and $\Gamma_{(p)}$ 
provides a notion of ``recurrence'' of a different sort. 
These groups each have a finite-index normal subgroup $K$ 
which naturally contains a direct product $K^{(\times p)}$ 
of copies of itself, again as a finite-index normal subgroup 
(for $\mathfrak{G}$ we take $p=2$). 
We may therefore consider the descending sequence 
$\Gamma_i = K^{(\times p^i)}$, 
and note that a commutator in $\Gamma_{i+1}$ is a $p$-tuple 
of commutators in $\Gamma_i$, 
thereby translating commutator calculations at the top few levels 
down to all other levels. 

The best previous diameter bounds for quotients of branch groups 
make no use of their branch structure, 
or indeed anything other than the fact that they are 
transitive permutation groups. 
Babai and Seress obtained the following very general result. 

\begin{thm}[\cite{BabSer} Theorem 1.4] \label{BabSerThm}
Let $G$ be a transitive permutation group of degree $N$. Then: 
\begin{center}
$\diam(G) \leq \exp (C (\log (N))^3) \diam \big(\Alt (m(G))\big)$
\end{center}
where $C$ is an absolute constant, 
and $m (G)$ is the maximal degree 
of an alternating composition factor of $G$. 
\end{thm}

If $\Gamma$ is a regular branch group acting on the $k$-ary rooted tree, then \linebreak$G = \Gamma / \Stab_{\Gamma} (n)$ 
has a natural transitive action on the $n$th level of the tree, 
so Theorem \ref{BabSerThm} is applicable, 
with $N = k^n$ and $m (G) \leq k$. 
For $\Gamma = \mathfrak{G}$ or $\Gamma_{(p)}$, $\Gamma / \Stab_{\Gamma} (n)$ is a $p$-group of order $\Omega (p^{\Omega (p^n)})$ 
(see Lemma \ref{stabindexlem} and Corollary \ref{GSStabindexcoroll} below), 
Theorem \ref{BabSerThm} yields:
\begin{align*}
\diam (\Gamma / \Stab_{\Gamma} (n)) 
& = O \Big(\exp \big(O (n^3 \log(p)^3)\big)\Big) \\
& = O \Big(\exp \big(O_p (\log \log \lvert \Gamma :\Stab_{\Gamma} (n) \rvert^3)\big)\Big)
\end{align*}
(here we take $p=2$ for $\Gamma=\mathfrak{G}$). Theorem \ref{BabSerThm} makes use of some deep machinery, 
including the Classification of Finite Simple Groups. 
Therefore our Theorems \ref{stabdiamcorollGrig}-\ref{centralcorollGS} 
improve upon prior results in at least three ways: 
\begin{itemize}
\item[(i)] They improve the diameter bound qualitatively, 
from a quasipolynomial function of $\log \lvert G \rvert$ 
to a polynomial one. 

\item[(ii)] They give explicit (and small) estimates for implied constants. 

\item[(iii)] They have self-contained, elementary and constructive proofs, which could in principle be implemented in reasonable time on a computer. 

\end{itemize}

It is also noteworthy that Theorem \ref{BabSerThm} remains a key 
tool in the study of other transitive permutation groups, 
which are consequently not known to have diameter less than 
$\exp (C (\log (N))^3)$. 
The best known bound in this direction is the following 
result of Helfgott and Seress. 

\begin{thm}[\cite{HelfSere}] \label{HelfSereThm}
Let $G$ be as in Theorem \ref{BabSerThm}. Then: 
\begin{center}
$\diam(G) \leq \exp (C (\log (N))^4 \log \log (N))$. 
\end{center}
\end{thm}

By Theorem \ref{BabSerThm}, Theorem \ref{HelfSereThm} 
can be immediately reduced to the case $G = \Alt (N)$. 
The proof in this special case also uses Theorem \ref{BabSerThm} 
to facilitate an important induction. 
It is a longstanding conjecture that $\Sym (N)$ and $\Alt (N)$ 
in reality have diameter polynomial in $N$ 
(that is, polylogarithmic in their order, 
like the groups studied in the present paper). 

A permutation group which provides a fascinating 
example intermediate between $\Sym (N)$ and 
the groups considered in 
Theorems \ref{stabdiamcorollGrig} and \ref{stabcorollGS}, 
is the Sylow $p$-subgroup 
$W_n = \Syl_p (\Sym (p^n))$ of $\Sym (p^n)$. 
$W_n$ was studied by Kaloujnine \cite{Kalou}, 
who showed that it is isomorphic to the 
$n$-fold iterated regular wreath product 
$C_p \wr C_p \wr \cdots \wr C_p$, 
which acts naturally on the (first $n$ levels of the) 
$p$-ary rooted tree. 
Although the inverse limit $\Gamma = \varprojlim_n W_n$ 
of the $W_n$ is a regular branch pro-$p$ group, 
it appears to be resistant to the Solovay-Kitaev procedure 
(for instance, $\Gamma$ is not finitely generated 
as a topological group; by contrast, 
every version of the profinite Solovay-Kitaev procedure 
known to the author proves finite generation as a byproduct). 
It is however likely that the methods of this paper are applicable 
to other residually nilpotent branch groups. 
This should be a topic of further study. 

The paper is structured as follows. 
Section \ref{prelimsect} is devoted to preliminary material: 
in subsection \ref{commsubsect} we lay out the consequences of 
the standard commutator identities upon which the profinite 
Solovay-Kitaev procedure is based, 
and illustrate, by means of an example, 
their relevance to diameter bounds. 
In subsection \ref{branchsection} we recall some basic material on 
group actions on regular rooted trees and the class of (regular) branch groups, define the Grigorchuk group 
$\mathfrak{G}$ and the Gupta-Sidki $p$-groups $\Gamma_{(p)}$ 
and give some basic properties. 
In Section \ref{mainthmGrigsect} we prove Theorem \ref{mainthmGrig}, 
and deduce Theorem \ref{stabdiamcorollGrig}. 
These proofs will be based in part upon prior results on the 
structure of the lower central series of $\mathfrak{G}$ 
(taken from \cite{Bar,BarGri}), which we state there. 
In Section \ref{GSproofsection} we prove Theorem \ref{mainthmGS} 
and deduce Theorem \ref{stabcorollGS}. 
In Section \ref{GS3sect} we recall some results from \cite{Bar} 
on the structure of the lower central series of $\Gamma_{(3)}$ 
and using them, deduce Theorem \ref{centralcorollGS} from 
Theorem \ref{mainthmGS}. 
In Section \ref{gapsect} we comment on implications of our diameter bounds 
for spectral gap and mixing times of random walks on Cayley graphs. 
Finally in Section \ref{growthsect} we make some remarks 
on the relationship between the growth of an infinite group 
and the diameters of its finite quotients. 

\section{Preliminaries} \label{prelimsect}

\subsection{Commutator Relations} \label{commsubsect}

The proofs of Theorems \ref{mainthmGrig} and \ref{mainthmGS} 
(from which our other results are deduced) 
will be via a variant of the 
profinite Solovay-Kitaev Procedure developed in \cite{Brad}, 
to which we refer the reader for further background. 
Given a group $G$ and a descending sequence of finite-index normal subgroups, 
the Procedure provides an approach to proving upper bounds on the diameters of the corresponding 
finite quotient groups, and relies on two key ingredients, 
both concerning the behaviour of commutator words within the group. 
The first ingredient encapsulates the intuitive notion that, 
given two specified elements $g,h \in G$ and two 
``good approximations'' $\tilde{g},\tilde{h} \in G$, 
the commutator $[\tilde{g},\tilde{h}]$ of the approximations is a good
approximation to the commutator $[g,h]$ of the original elements. 
What is perhaps not so intuitive is that in many situations, 
$[\tilde{g},\tilde{h}]$ approximates $[g,h]$ much more closely than 
$\tilde{g}$ and $\tilde{h}$ did $g$ and $h$. 
This idea is made precise in the following Lemma, 
which is an immediate consequence of standard commutator identities. 

\begin{lem} \label{SKlem1}
Let $G_1 , G_2 , H_1 , H_2 \vartriangleleft G$, with $G_1 \geq G_2$, $H_1 \geq H_2$. 
For all $g_i \in G_i$, $h_i \in H_i$ ($i=1,2$), 
\begin{center}
$[g_1 g_2 , h_1 h_2] \equiv [g_1 , h_1] \mod [G_1,H_2] [G_2,H_1]$. 
\end{center}
\end{lem}

\begin{proof}
We compute directly: 
\begin{center}
$[g_1 g_2,h_1 h_2] = [g_1,h_2][g_1,h_1][[g_1,h_1],h_2][[g_1,h_1 h_2],g_2][g_2,h_1 h_2]$
\end{center}
and all terms other than $[g_1,h_1]$ lie in $[G_1,H_2] [G_2,H_1]$. 
\end{proof}

A particularly useful special case of this lemma is the following. 

\begin{coroll} \label{SKlem2}
Let $(K_i)_{n=1} ^{\infty}$ be a descending sequence of normal subgroups of $G$. 
Suppose that, for all $m,n \in \mathbb{N}$, $[K_m,K_n] \subseteq K_{m+n}$. 
Let $m_1,m_2,n_1,n_2 \in \mathbb{N}$, with $m_1 \leq m_2$, $n_1 \leq n_2$, 
and let $g_i \in K_{m_i}$, $h_i \in K_{n_i}$ ($i=1,2$). Then:
\begin{center}
$[g_1 g_2 , h_1 h_2] \equiv [g_1 , h_1] \mod K_{\min (m_1 + n_2 , m_2 + n_1)}$. 
\end{center}
\end{coroll}

\begin{proof}
Set $G_i = K_{m_i}$ and $H_i = K_{n_i}$ in Lemma \ref{SKlem1}, for $i=1,2$. 
\end{proof}

We will use Corollary \ref{SKlem2} in the proof of Theorem \ref{mainthmGrig}. 
It is a classical fact that the lower central series satisfies the required hypothesis. 

\begin{lem} \label{lcentcommlem}
Let $G$ be any group. Then for all $m,n \in \mathbb{N}$, 
\begin{center}
$[\gamma_n (G),\gamma_m (G)] \subseteq \gamma_{m+n} (G)$. 
\end{center}
\end{lem}

Lemma \ref{lcentcommlem} will be used extensively and without further comment in the sequel, and in particular throughout Section \ref{mainthmGrigsect}. 

What is the relevance of the preceding discussion to diameters? 
As we intimated in the introduction, 
given a group $\Gamma$ and a descending sequence $(\Gamma_i)_i$, 
the Solovay-Kitaev procedure requires the existence 
of an approximation to elements of a deeper term $\Gamma_j$ in the sequence, 
by commutators $[g,h]$ of elements $g,h$ lying in a higher term $\Gamma_i$. 
We may assume by induction that we have approximations $\tilde{g},\tilde{h}$ 
to $g,h$ up to an error lying in the intermediate term $\Gamma_k$, 
where $\tilde{g}$ and $\tilde{h}$ are short words in a generating set. 
Substituting $\tilde{g}$ and $\tilde{h}$ into our commutator expression, 
we obtain approximations $[\tilde{g},\tilde{h}]$ to elements of $\Gamma_j$ 
by short words. Lemma \ref{SKlem1} 
gives us some control over the fidelity of these approximations. 
Let us illustrate this with an example. 

\begin{ex} \label{SKExample}
\normalfont 
Let $\Gamma$ be a group, let $\Gamma \geq \Gamma_1 \geq \Gamma_2 \geq \Gamma_3$ be finite-index normal subgroups of $\Gamma$, and suppose: 
\begin{equation} \label{shrinkeqnex}
[\Gamma_1 , \Gamma_2] \leq \Gamma_3\text{.}
\end{equation}
Second, suppose that for any $k \in \Gamma_2$, there exist $g,h \in \Gamma_1$ 
such that: 
\begin{equation} \label{approxeqnex}
[g,h]k^{-1} \in \Gamma_3\text{.} 
\end{equation}
Now let $S \subseteq \Gamma$, 
and suppose the image of $S$ in $\Gamma/\Gamma_3$ is a generating set. 
Let $d = \diam (\Gamma/\Gamma_2,S)$, so that: 
\begin{equation} \label{hypeqnex}
\Gamma = \Gamma_2 B_S (d) 
\end{equation}
(here we abuse notation slightly, denoting also by ``$S$'' the image of $S$ 
modulo any of the $\Gamma_i$). 

We wish to bound $\diam (\Gamma/\Gamma_3,S)$. 
That is, given $k \in \Gamma$ we seek a short word $w$ in $S$ 
such that $k w^{-1} \in \Gamma_3$. 
First suppose that $k \in \Gamma_2$. 
Let $g,h \in K_1$ be as in (\ref{approxeqnex}). 
By (\ref{hypeqnex}), there exist $\tilde{g} , \tilde{h} \in B_S (d)$ such that 
$g^{-1} \tilde{g} , h^{-1} \tilde{h} \in \Gamma_2$. 
Setting $G_1 = H_1 = \Gamma_1$, $G_2 = H_2 = \Gamma_2$, 
$g_1 = g$, $h_1 = h$, $g_2 = g^{-1} \tilde{g}$, $h_2 = h^{-1} \tilde{h}$ 
in Lemma \ref{SKlem1}, and applying (\ref{shrinkeqnex}), 
\begin{center}
$[\tilde{g},\tilde{h}] \equiv [g,h] \equiv k \mod \Gamma_3$. 
\end{center}
Since $[\tilde{g},\tilde{h}] \in B_S (4d)$, we have 
$\Gamma_2 \subseteq \Gamma_3 B_S (4d)$. Applying (\ref{hypeqnex}) again, 
$\Gamma = \Gamma_3 B_S (5d)$, so $\diam(\Gamma/\Gamma_3,S) \leq 5d$. 
\end{ex}

Although the hypotheses (\ref{shrinkeqnex}) and (\ref{approxeqnex}) 
may seem somewhat artificial in the context of this abstract example, 
in practice many groups satisfy these conditions, or variants thereof. 
Modifications to the method of Example \ref{SKExample} 
are however sometimes necessary, or desirable. 
For instance the elements $g$ and $h$ in Example \ref{SKExample} 
were taken to lie at the same level of the subgroup chain, 
whereas for many groups, including the branch groups studied in this paper, 
good commutator expressions will involve elements lying at different levels. 

\subsection{Groups Acting on Regular Rooted Trees} \label{branchsection}

Let us set up some basic notation and definitions concerning the class of groups to be studied in the sequel. 
All of the material here (and much, much more besides) is covered in \cite{BaGrSu} and \cite{dlHa} Chapter VIII. 

Let $\mathcal{A}$ be a finite set, 
and let $\mathcal{A} ^*$ be the set of formal positive words on the alphabet $\mathcal{A}$. 
We may partially order $\mathcal{A}^*$ via the \emph{prefix relation} $\leq$, 
where for $u,v \in \mathcal{A}^*$, $u \leq v$ iff there exists $w \in \mathcal{A}^*$ 
such that $v=uw$. 

Geometrically, we may regard $\mathcal{A}^*$ 
as the set of vertices of a regular rooted tree $\mathcal{T}_{\mathcal{A}}$: 
the root vertex is identified with the empty word, 
and every vertex $v$ is joined by an edge to its $\lvert \mathcal{A} \rvert$ \emph{children} 
$v a$, $a \in \mathcal{A}$. 
Under this identification, the set $\mathcal{A}^{n}$ of words of length $n$ 
is precisely the sphere of radius $n$ in $\mathcal{T}_{\mathcal{A}}$ about the root vertex, 
known as the \emph{$n$th level set}. 

The \emph{automorphism group} $\Aut (\mathcal{T}_{\mathcal{A}})$ 
is the set of permutations of $\mathcal{A}^*$ 
preserving the prefix relation. 
Geometrically this just the group of graph automorphisms of the tree $\mathcal{T}_{\mathcal{A}}$. 
For the valence of every vertex of $\mathcal{T}_{\mathcal{A}}$ is $\lvert \mathcal{A} \rvert + 1$ 
except for the root vertex (which has valence $\lvert \mathcal{A} \rvert$), 
so every automorphism of $\mathcal{T}_{\mathcal{A}}$ fixes the root vertex, 
and hence preserves the level sets. 
The kernel of the action of $\Aut (\mathcal{T}_{\mathcal{A}})$ on the $n$th level set $\mathcal{A}^{n}$ 
will be called the \emph{$n$th level stabiliser} and denoted $\Stab (n)$; 
it is naturally isomorphic to $\Aut (\mathcal{T}_{\mathcal{A}})^{(\times \lvert \mathcal{A} \rvert^n)}$. 
If $\Gamma \leq \Aut (\mathcal{T}_{\mathcal{A}})$ we write $\Stab_{\Gamma} (n)$ for $\Gamma \cap \Stab (n)$, 
though in general we cannot say more about the structure of $\Stab_{\Gamma} (n)$ than that 
$\Stab_{\Gamma} (n)$ is isomorphic to a subgroup of 
$\Aut (\mathcal{T}_{\mathcal{A}})^{(\times \lvert \mathcal{A} \rvert^n)}$. 

For any $\phi \in \Aut (\mathcal{T}_{\mathcal{A}})$, 
there exists a unique $\sigma_{\phi} \in \Sym (\mathcal{A})$ 
such that for any $x \in \mathcal{A}$, there exists a unique $\phi_x \in \Aut (\mathcal{T}_{\mathcal{A}})$ such that: 
\begin{center}
$\phi (x w) = \sigma_{\phi} (x) \phi_x (w)$, for all $w \in \mathcal{A}^*$. 
\end{center}
The induced map $\psi : \phi \mapsto (\phi_x)_{x \in \mathcal{A}} \cdot \sigma_{\phi}$ 
gives an isomorphism \linebreak$\Aut (\mathcal{T}_{\mathcal{A}}) 
\rightarrow \Aut (\mathcal{T}_{\mathcal{A}}) \wr \Sym (\mathcal{A})$. 
Note that the level stabilisers may be described recursively by 
$\Stab (0) = \Aut (\mathcal{T}_{\mathcal{A}})$ and 
$\Stab (n+1) = \psi^{-1} (\Stab (n)^{(\times \lvert \mathcal{A} \rvert)})$ for $n \in \mathbb{N}$. 

Of particular interest among the subgroups of $\Aut (\mathcal{T}_{\mathcal{A}})$ 
are those whose action on $\mathcal{A}^*$ is \emph{branch}. 
Our characterization of such groups is based on that appearing in \cite{Bar}. 

\begin{defn} \label{branchdefn}
Let $\Gamma \leq \Aut (\mathcal{T}_{\mathcal{A}})$. 
$\Gamma$ is \emph{(regular) branch} if: 
\begin{itemize}
\item[(i)] The action of $\Gamma$ on $\mathcal{A}$ is transitive; 
\item[(ii)] $\psi (\Stab_{\Gamma} (1)) \leq \Gamma^{(\times \lvert \mathcal{A} \rvert)}$; 
\item[(iii)] $\Gamma$ has a finite-index subgroup $K$ such that $K^{(\times \lvert \mathcal{A} \rvert)} \leq \psi(K)$. 
\end{itemize}
We will simply say that a group $\Gamma$ \emph{branches over $K$} when the alphabet $\mathcal{A}$ 
and the action of $\Gamma$ on $\mathcal{A}^*$ is clear. 
\end{defn}

For the sake of uncluttered notation 
we will allow ourselves to supress the map $\psi$ from expressions and identify subgroups of $\Gamma$ 
with their image under $\psi$, 
so for instance we may (abuse notation somewhat and) speak of 
$K^{(\times \lvert \mathcal{A} \rvert)}$ as a subgroup of $K$; 
$\Stab_{\Gamma} (n)$ as a subgroup of $\Gamma ^{(\times \lvert \mathcal{A} \rvert^n)}$ and so on. 

In truth, branch groups form a much broader class than those groups covered by Definition \ref{branchdefn}
(see for instance \cite{BaGrSu}), but this more restricted setting will be most convenient for our purposes. 

We now define the specific branch groups which are the subject of Theorems \ref{stabdiamcorollGrig}-\ref{centralcorollGS}. 

\subsubsection{The Grigorchuk $2$-Group} \label{Grigdefnsubsubsect}

Let $\mathcal{A} = \lbrace 0,1 \rbrace$ and write 
$\mathcal{T}_{\mathcal{A}} = \mathcal{T}_2$. 
The \emph{Grigorchuk $2$-group} 
(sometimes known as the \emph{first Grigorchuk group}) 
is the subgroup $\mathfrak{G}$ of $\Aut (\mathcal{T}_2)$ 
generated by the four automorphisms $a,b,c,d$, defined by: 
\begin{center}
$a (0 w) = 1 w$; $a (1 w) = 0 w$; \\
$b (0 w) = 0 a (w)$; $b (1 w) = 1 c(w)$; \\
$c (0 w) = 0 a (w)$; $c (1 w) = 1 d(w)$; \\
$d (0 w) = 0 w$; $d (1 w) = 1 b(w)$. 
\end{center}
In other words, $a$ swaps the subtrees rooted at $0$ and $1$, 
while $b,c,d \in \Stab_{\mathfrak{G}} (1)$ are defined recursively via: 
\begin{center}
$b = (a,c), c = (a,d), d = (1,b)$. 
\end{center}
An easy induction on the levels shows that: 
\begin{equation} \label{Grigrelns}
a^2 = b^2 = c^2 = d^2 = 1, bc=cb=d, cd=dc=b, bd=db=c\text{.}
\end{equation}
Let $x = abab \in \mathfrak{G}$ and let $K = \langle x \rangle^{\mathfrak{G}} \vartriangleleft \mathfrak{G}$. 
We will require the following basic facts in the sequel. 

\begin{propn}[\cite{dlHa} Chapter VIII] \label{factorsareC4}
Let $\mathfrak{G}$, $K$, $x$ be as above. 
\begin{itemize}
\item[(i)] $\mathfrak{G}$ branches over $K$;
\item[(ii)] $\mathfrak{G} / K \cong D_8 \times C_2$; 
\item[(iii)] $K / K^{(\times 2)} \cong C_4$, generated by $x$. 
\end{itemize}

\end{propn}

\begin{lem}
$K^{(\times 2^m)} \vartriangleleft \mathfrak{G}$ for all $m \geq 0$. 
\end{lem}

\begin{proof}
We proceed by induction on $m$, the base case $m=0$ holding by definition. 

For $m \geq 1$, we verify that $K^{(\times 2^m)}$ is preserved under conjugation by the generators $a,b,c,d$. 
This is the case for $a$, which simply permutes the factors in the direct product. 

$b,c$ and $d$ preserve the decomposition $K^{(\times 2^m)} = (K^{(\times 2^{m-1})})^{(\times 2)}$, 
and act on each $K^{(\times 2^{m-1})}$-factor as $a,b,c$ or $d$. The result follows by induction. 
\end{proof}

\begin{lem} \label{Kinstablem}
For all $m \geq 0$, $K^{(\times 2^m)} \leq \Stab_{\mathfrak{G}} (m+1)$. 
\end{lem}

\begin{proof}
Note that $x = abab \in \Stab_{\mathfrak{G}} (1)$. The result is now immediate from Proposition \ref{factorsareC4} (iii). 
\end{proof}

\begin{lem} \label{stabindexlem}
For all $n \geq 1$, 
\begin{center}
$\lvert \mathfrak{G} : \Stab_{\mathfrak{G}} (n) \rvert 
\geq 2^{2^{n-1}+1}$. 
\end{center}
\end{lem}

\begin{proof}
It suffices to prove that 
$\lvert \Stab_{\mathfrak{G}} (m) : \Stab_{\mathfrak{G}} (m+1) \rvert 
\geq 2^{2^{m-1}}$ for all $m \in \mathbb{N}$. 
$x=abab \equiv (a,a) \mod \Stab_{\Aut (\mathcal{T}_2)} (2)$, 
so as $\underline{\epsilon}$ ranges over $\lbrace 0,1 \rbrace^{2^m}$, 
the elements $(x^{\epsilon_i})_{i=1} ^{2^m}$ are all distinct 
modulo $\Stab_{\mathfrak{G}} (m+1)$, 
and lie in $\Stab_{\mathfrak{G}} (m)$ by Lemma \ref{Kinstablem}. 
\end{proof}

\subsubsection{The Gupta-Sidki $p$-Groups}

Fix an odd prime $p$, 
let $\mathcal{A} = \lbrace 0,1, \ldots , p-1 \rbrace$ 
and write $\mathcal{T}_{\mathcal{A}} = \mathcal{T}_p$. 
The \emph{Gupta-Sidki $p$-group} is the subgroup $\Gamma_{(p)}$ of $\Aut (\mathcal{T}_p)$ 
generated by the two automorphisms $a,b$, where $a$ is defined by: 
\begin{center}
$a (i w) = (i+1) w$ for $0 \leq i \leq p-2 $; $a ((p-1) w) = 0 w$
\end{center}
(so that $a$ cyclically permutes the level-$1$ subtrees) and $b \in \Stab_{\Aut (\mathcal{T}_p)} (1)$ is defined recursively via: 
\begin{center}
$b = (a,a^{-1},1,\ldots,1,b)$ 
\end{center}
so that both $a$ and $b$ have order $p$. 
Note that a reduced word in $a$ and $b$ corresponds to an element of 
$\Stab_{\Gamma_{(p)}} (1)$ iff the number of occurences of $a$ 
(counted with signs) is congruent to $0$ modulo $p$. 
Thus $\Stab_{\Gamma_{(p)}} (1) = \langle b \rangle^{\Gamma_{(p)}}$. 

Let $K = [\Gamma_{(p)},\Gamma_{(p)}]$ be the derived subgroup of $\Gamma_{(p)}$. Note that \linebreak$K \leq \Stab_{\Gamma_{(p)}} (1)$. 
Let $x_1 = [a,b] \in K$, and for $1 \leq i \leq p-2$, define $x_{i+1}=[a,x_i] \in K$. 
The following computations were made in \cite{Garr}, 
generalising results from \cite{Sidki} for $\Gamma_{(3)}$. 

\begin{propn}[\cite{Garr} Proposition 2.2] \label{Garridoprop}
Let $B = \Stab_{\Gamma_{(p)}} (1)$. Then: 
\begin{itemize}
\item[(i)] $\Gamma_{(p)} / K \cong C_p \times C_p$, with basis $Ka,Kb$; 

\item[(ii)] $B^{\prime} = K^{(\times p)}$; 

\item[(iii)] $\Gamma_{(p)} / B^{\prime} \cong C_p \wr C_p$. 

\end{itemize}
\end{propn}

From this we have an explicit description 
of the branch structure of $\Gamma_{(p)}$. 

\begin{coroll} \label{GSbranchprop}
\begin{itemize}
\item[(i)] $\Gamma_{(p)}$ branches over $K$; 

\item[(ii)] $K / K^{(\times p)} \cong C_p ^{(\times (p-1))}$, with basis $x_1 , \ldots , x_{p-1}$. 

\end{itemize}
\end{coroll}

\begin{proof}
\begin{itemize}
\item[(i)] From Proposition \ref{Garridoprop} (ii), 
\begin{center}
$K^{(\times p)} = B^{\prime} \leq \Gamma_{(p)} ^{\prime} = K$. 
\end{center}

\item[(ii)] We defer this to Section \ref{GSproofsection}, 
where we introduce additional notation which will be 
convenient to use in the proof. 
\end{itemize}
\end{proof}

The following observation will be key to the deduction 
of Theorem \ref{stabcorollGS} from Theorem \ref{mainthmGS}. 

\begin{lem} \label{GSKStabincl}
For all $m \in \mathbb{N}$, 
$K^{(\times p^m)} \leq \Stab_{\Gamma_{(p)}} (m+1)$. 
\end{lem}

\begin{proof}
It clearly suffices to check that $K \leq \Stab_{\Gamma_{(p)}}(1)$. 
This is so because $K = [\Gamma_{(p)},\Gamma_{(p)}]$ 
and $\Gamma_{(p)} / \Stab_{\Gamma_{(p)}}(1)$ is abelian. 
\end{proof}

\section{Proofs for Grigorchuk's Group} \label{mainthmGrigsect}

In this Section we prove Theorem \ref{mainthmGrig} 
and deduce Theorem \ref{stabdiamcorollGrig}. 
Before embarking on the proof of our diameter bounds, 
we marshal some facts about the lower central series of $\mathfrak{G}$. 

Recall that for any group $G$, the degree $\deg (g)$ of $g \in G$ is given by 
\linebreak$g \in \gamma_{\deg(g)} (G) \setminus \gamma_{\deg(g)+1} (G)$, with the convention that all $g \in \bigcap_{n=1} ^{\infty} \gamma_n (G)$ have degree $\infty$ (though in $\mathfrak{G}$, the latter situation never arises for $g \neq 1$: $\mathfrak{G}$ is a residually finite $2$-group, so is residually $2$-finite, and in particular, $\bigcap_{n=1} ^{\infty} \gamma_n (\mathfrak{G}) = \lbrace 1 \rbrace$). 

For any $g \in K$, write $\mathbf{0}(g) = (g,1), \mathbf{1}(g) = (g,g^{-1}) \in K^{(\times 2)}$. 

\begin{thm}[\cite{Bar,BarGri,Roz}] \label{degreethmGrig}
Let $X_1 , \ldots , X_n \in \lbrace \mathbf{0},\mathbf{1} \rbrace$. Then: 
\begin{equation*}
\begin{aligned}
\deg X_1 \cdots X_n (x) & = 1 + \sum_{i=1} ^n X_i 2^{i-1} + 2^n\text{;} \\
\deg X_1 \cdots X_n (x^2) & = 1 + \sum_{i=1} ^n X_i 2^{i-1} + 2^{n+1}\text{.} 
\end{aligned}
\end{equation*}
\end{thm}

\begin{thm}[\cite{Bar,BarGri,Roz}] \label{centralindexGrig}
For all $m \geq 2$, 
\begin{center}
$\lvert \mathfrak{G} : \gamma_{2^m + 1} (\mathfrak{G}) \rvert 
= 2 ^{2^{m-1} 3 + 2}$. 
\end{center}
\end{thm}

To avoid cluttered notation, in this section we may write 
$\gamma_n = \gamma_n (\mathfrak{G})$. 
The specific consequences of Theorem \ref{degreethmGrig} 
to be used in the proof of Theorem \ref{mainthmGrig} are as follows. 

\begin{coroll} \label{degreecorollKs}
$\gamma_{2^m + 2^{m-1} + 1} \leq K^{(\times 2^m)} \leq \gamma_{2^m + 1}$. 
\end{coroll}

To be more precise, we have the following estimates. 

\begin{coroll} \label{degcompcoroll}
For all $\underline{\delta} \in \big(\lbrace 0,1 \rbrace^{2^{n-2}}\big) \setminus \lbrace \underline{0} \rbrace$, 
\begin{itemize}
\item[(i)] $2^n + 1 \leq \deg \big( \big((x,1,1,1)^{\delta_i}\big)_{i=1} ^{2^{n-2}} \big) \leq 2^n+2^{n-2}$; 

\item[(ii)] $2^n +2^{n-2}+ 1 \leq \deg \big( \big((x,1,x,1)^{\delta_i}\big)_{i=1} ^{2^{n-2}} \big) \leq 2^n +2^{n-1}$;

\item[(iii)] $2^n +2^{n-1}+ 1 \leq \deg \big( \big((x,x,1,1)^{\delta_i}\big)_{i=1} ^{2^{n-2}} \big) 
\leq 2^n +3 \cdot 2^{n-2}$;

\item[(iv)] $2^n +3 \cdot 2^{n-2}+ 1 \leq \deg \big( \big((x,x,x,x)^{\delta_i}\big)_{i=1} ^{2^{n-2}} \big) 
\leq 2^{n+1}$;

\item[(v)] Specifically, 
$\deg \big( \big((x,x,x,x)\big)_{i=1} ^{2^{n-2}} \big) = 2^{n+1}$; 

\item[(vi)] $2^n + 1 \leq \deg \big( \big((x^2,1)^{\delta_i}\big)_{i=1} ^{2^{n-2}} \big) \leq 2^n+2^{n-2}$; 

\item[(vii)] $2^n+2^{n-2} + 1 \leq \deg \big( \big((x^2,x^2)^{\delta_i}\big)_{i=1} ^{2^{n-2}} \big) \leq 2^n +2^{n-1}$; 

\item[(viii)] Specifically, 
$\deg \big( \big((x^2,x^2)\big)_{i=1} ^{2^{n-2}} \big) = 2^n +2^{n-1}$. 

\end{itemize}
\end{coroll}

We now produce our commutator approximations to elements lying sufficiently 
deep in $(\gamma_n)_n$. 
The following identities are verified by direct computation. 

\begin{lem} \label{commscompGriglem}
\begin{itemize}
\item[(i)] $\big[ x,(x,1) \big] 
= \big( x^{-1},1,1,1 \big)$. 

\item[(ii)] $\big[ x,(x,x) \big] 
= \big( x^{-1},1,1,(1,x^{-1})x \big)$. 

\item[(iii)] $\big[ x^2,(x,1) \big] 
= \big( x^{-1},x,1,1 \big)$. 

\item[(iv)] $\big[ x^2,(x,x^{-1}) \big] 
= \big( x^{-1},x,(x^{-1},1)x^{-1},(1,x^{-1})x \big)$. 

\end{itemize}
\end{lem}

\begin{proof}
We prove (i) and leave the verifications of the other identities, 
which are similar, as an exercise. 
First note that $x = abab = (ca,ac)$. Thus: 
\begin{center}
$\big[ x,(x,1) \big] = \big( [ca,x] , 1 \big)$
\end{center}
Now, using (\ref{Grigrelns}) we have: 
\begin{align*}
[ca,x] & = ac(ac,ca)ca(ca,ac) \\
 & = a(ca,bad)a(ca,ac) \\
  & = (bad,ca)(ca,ac) \\
   & = (baba,1) \\
      & = (x^{-1},1) 
\end{align*}
as desired. 
\end{proof}

Lemma \ref{commscompGriglem} yields good approximations to 
elements of $K^{(\times 4)}$, modulo $K^{(\times 8)}$, 
by commutators of elements of $K$ and $K^{(\times 2)}$. 
For deeper subgroups $K^{(\times 2^n)}$, 
we express elements as vectors of elements in $K^{(\times 4)}$, 
and produce a commutator approximation by applying 
Lemma \ref{commscompGriglem} to each term of the vector. 
The expressions we obtain can be related to the lower central series 
by using Theorem \ref{degreethmGrig} and Corollaries 
\ref{degreecorollKs} and \ref{degcompcoroll} to estimate the 
degrees of the elements occuring. 
In summary we have the following Proposition. 

\begin{propn} \label{commspropnGrig}
Let $\mathfrak{c} : \mathfrak{G} \times \mathfrak{G} \rightarrow \mathfrak{G}$ be given by $\mathfrak{c}(g,h)=[g,h]$. 
Let $m \geq 2$. Then: 
\begin{itemize}
\item[(i)] The restriction of $\mathfrak{c}$ to $\gamma_{2^{m-1}} \times \gamma_{2^{m-1}}$ descends to a well-defined map:
\begin{center}
$\bar{\mathfrak{c}}_m : (\gamma_{2^{m-1}} / \gamma_{2^m + 1})^{(\times 2)} 
\rightarrow \mathfrak{G} / \gamma_{2^m + 2^{m-1} + 1}$
\end{center}
whose image contains $K^{(\times 2^m)} / \gamma_{2^m + 2^{m-1} + 1}$. 
\item[(ii)] The restriction of $\mathfrak{c}$ to 
$\gamma_{2^{m-1} + 2^{m-2}} \times \gamma_{2^{m-1} + 2^{m-2}}$ 
descends to a well-defined map:
\begin{center}
$\bar{\bar{\mathfrak{c}}}_m : (\gamma_{2^{m-1} + 2^{m-2}} / \gamma_{2^m + 2^{m-1} + 1})^{(\times 2)}
\rightarrow \mathfrak{G} / \gamma_{2^{m+1} + 1}$
\end{center}
whose image contains $\gamma_{2^m + 2^{m-1} + 1} / \gamma_{2^{m+1} + 1}$. 
\end{itemize}
\end{propn}

\begin{proof}
The well-definedness of $\bar{\mathfrak{c}}_m$ 
and $\bar{\bar{\mathfrak{c}}}_m$ is an immediate consequence of 
Corollary \ref{SKlem2}. It therefore suffices to check the images 
of the maps contain the specified subgroups. 

For (i), note that by Corollary \ref{degcompcoroll}, 
every element of $K^{(\times 2^m)} / \gamma_{2^m + 2^{m-1} + 1}$ 
is represented by 
$\big( (x^{\delta_i + \epsilon_i},1,x^{\epsilon_i},1) \big)_{i=1} ^{2^{m-2}}$ 
for some $\delta_i , \epsilon_i \in \lbrace 0,1 \rbrace$. 
By Lemma \ref{commscompGriglem} (i) and (ii), 
\begin{equation*}
\begin{aligned}
\big[ x,(x,1) \big] \cdot \big( x,1,1,1 \big)^{-1} & = \big( x^2,1,1,1 \big)^{-1}  \\
\big[ x,(x,x) \big] \cdot \big( x,1,x,1 \big)^{-1} & \equiv 
\big( 1,1,x,x \big) \big( x^2,1,x^2,1 \big)^{-1} \mod K^{(\times 2^3)}\text{.}
\end{aligned} 
\end{equation*}
By Corollary \ref{degcompcoroll} (iii) and (iv), for any $(\beta_i)_{i=1} ^{2^{m-2}} \in \lbrace 0,1 \rbrace^{2^{m-2}}$, 
\begin{center}
$\big( (1,1,x,x)^{\beta_i} \big)_{i=1} ^{2^{m-2}} \in \gamma_{2^m + 2^{m-1} + 1}$, 
\end{center}
by Corollary \ref{degcompcoroll} (vi), for any $(\beta_i)_{i=1} ^{2^{m-2}} \in \lbrace 0,1 \rbrace^{2^{m-2}}$, 
\begin{center}
$\big( (x^2,1,1,1)^{\beta_i} \big)_{i=1} ^{2^{m-2}}, 
\big( (x^2,1,x^2,1)^{\beta_i} \big)_{i=1} ^{2^{m-2}} 
\in \gamma_{2^m + 2^{m-1} + 1}$, 

\end{center}
and from Corollary \ref{degreecorollKs}, $K^{(\times 2^{m+1})} \subseteq \gamma_{2^{m+1}+1}$. 
Hence, for any $(\delta_i)_{i=1} ^{2^{m-2}}$, \linebreak$(\epsilon_i)_{i=1} ^{2^{m-2}} \in \lbrace 0,1 \rbrace^{2^{m-2}}$, 
\begin{equation*} 
\begin{aligned}
\big( (x^{\delta_i + \epsilon_i},1,x^{\epsilon_i},1) \big)_{i=1} ^{2^{m-2}} 
 \equiv & \big[ (x)_{i=1} ^{2^{m-2}} , ((x,1)^{\delta_i})_{i=1} ^{2^{m-2}} \big] \\
& \cdot \big[ (x)_{i=1} ^{2^{m-2}} , ((x,x)^{\epsilon_i})_{i=1} ^{2^{m-2}} \big] \mod \gamma_{2^m + 2^{m-1} + 1}\text{.}
\end{aligned}
\end{equation*} 
Using Corollary \ref{degcompcoroll} to estimate the degrees of $(x)_{i=1} ^{2^{m-2}}$, 
$((x,1)^{\delta_i})_{i=1} ^{2^{m-2}}$ and $((x,x)^{\epsilon_i})_{i=1} ^{2^{m-2}}$, 
and by the standard identity $[a,bc]=[a,c][a,b][[a,b],c]$, we deduce: 
\begin{equation} \label{Grigeqn1}
\big( (x^{\delta_i + \epsilon_i},1,x^{\epsilon_i},1) \big)_{i=1} ^{2^{m-2}}
\equiv \big[ (x)_{i=1} ^{2^{m-2}} , ((x^{\delta_i + \epsilon_i},x^{\epsilon_i}))_{i=1} ^{2^{m-2}} \big] 
\mod \gamma_{2^m + 2^{m-1} + 1}\text{.}
\end{equation}
For (ii), we see similarly by Corollary \ref{degcompcoroll} 
that every element of $\gamma_{2^m + 2^{m-1} + 1} / \gamma_{2^{m+1} + 1}$ 
is represented by 
$\big( (x^{\delta_i + \epsilon_i},x^{\delta_i + \epsilon_i},x^{\epsilon_i},x^{\epsilon_i}) \big)_{i=1} ^{2^{m-2}}$ 
for some $\delta_i , \epsilon_i \in \lbrace 0,1 \rbrace$. 
By Lemma \ref{commscompGriglem} (iii) and (iv), 
\begin{equation*}
\begin{aligned}
\big[ x^2,(x,1) \big] \cdot \big( x,x,1,1 \big)^{-1} & = \big( x^2,1,1,1 \big)^{-1} \\
\big[ x^2,(x,x^{-1}) \big] \cdot \big( x,x,x,x \big)^{-1} & \equiv \big( x^2,1,x^2,1 \big)^{-1} \mod K^{(\times 2^3)}\text{.}
\end{aligned}
\end{equation*}
Using Corollary \ref{degreecorollKs} 
and Corollary  \ref{degcompcoroll} (vi) 
once again as in (i), we have that for any $(\delta_i)_{i=1} ^{2^{m-2}} , (\epsilon_i)_{i=1} ^{2^{m-2}} \in \lbrace 0,1 \rbrace^{2^{m-2}}$, 
\begin{equation*}
\begin{aligned}
\big( (x^{\delta_i + \epsilon_i},x^{\delta_i + \epsilon_i},x^{\epsilon_i},x^{\epsilon_i}) \big)_{i=1} ^{2^{m-2}} 
\equiv & \big[ (x^2)_{i=1} ^{2^{m-2}} , ((x,1)^{\delta_i})_{i=1} ^{2^{m-2}} \big] \\ & \cdot
\big[ (x^2)_{i=1} ^{2^{m-2}} , ((x,x^{-1})^{\epsilon_i})_{i=1} ^{2^{m-2}} \big] \mod \gamma_{2^m + 2^{m-1} + 1}\text{.}
\end{aligned}
\end{equation*}
As before, we apply the commutator identity for products, 
using the estimate of the degrees of $(x^2)_{i=1} ^{2^{m-2}}$, 
$((x,1)^{\delta_i})_{i=1} ^{2^{m-2}}$ and $((x,x^{-1})^{\epsilon_i})_{i=1} ^{2^{m-2}}$ from Corollary \ref{degcompcoroll}, 
and deduce: 
\begin{equation} \label{Grigeqn4}
\big( (x^{\delta_i + \epsilon_i},x^{\delta_i + \epsilon_i},x^{\epsilon_i},x^{\epsilon_i}) \big)_{i=1} ^{2^{m-2}} 
\equiv \big[ (x^2)_{i=1} ^{2^{m-2}} , ((x^{\delta_i + \epsilon_i},x^{-\epsilon_i}))_{i=1} ^{2^{m-2}} \big] 
\mod \gamma_{2^{m+1} + 1}\text{.}
\end{equation}
\end{proof}

Using Proposition \ref{commspropnGrig} we may approximate 
any tuple consisting of $x$s and $1$s by a commutator. 
To express arbitrary tuples in 
$K^{(\times 2^{m-1})}/K^{(\times 2^{m})}$ we must also 
find approximations for tuples consisting of $x^2$s and $1$s. 
Here we will diverge slightly from our overall strategy of 
approximating elements by commutators, 
since it appears far more natural to express such tuples as 
squares. 

\begin{propn} \label{squaring}
The squaring map $\mathfrak{s} : g \mapsto g^2$ on $\mathfrak{G}$ induces a surjection: 
\begin{center}
$\overline{\mathfrak{s}}_m:K^{(\times 2^{m-1})}/\gamma_{2^m + 1} \rightarrow\gamma_{2^m + 1} /K^{(\times 2^m)}$
\end{center}
for all $m \geq 1$. 
\end{propn}

\begin{proof}
It is immediate from Lemma \ref{factorsareC4} and Corollaries \ref{degreecorollKs} and \ref{degcompcoroll} that 
\linebreak$K^{(\times 2^{m-1})}/\gamma_{2^m + 1}$ and 
$\gamma_{2^m + 1} / K^{(\times 2^m)}$ are elementary abelian $2$-groups, 
each element of the latter being represented by a vector: 
\begin{center}
$(x^{2 \delta_i})_{i=1} ^{2^{m-1}} = \mathfrak{s}((x^{\delta_i})_{i=1} ^{2^{m-1}})$,
\end{center}
as $\delta$ ranges over $\lbrace 0,1 \rbrace^{2^{m-1}}$. 
Since $(x^{\delta_i})_{i=1} ^{2^{m-1}} \in K^{(\times 2^{m-1})}$, 
$\mathfrak{s}$ induces a surjection $K^{(\times 2^{m-1})} \rightarrow\gamma_{2^m + 1} /K^{(\times 2^m)}$. 
\\ \\
Now let $a \in K^{(\times 2^{m-1})}$, $b \in \gamma_{2^m + 1}$. 
We have: 
\begin{center}
$\mathfrak{s} (ab) = \mathfrak{s} (a) [a,b^{-1}] \mathfrak{s} (b)$. 
\end{center}
But $\mathfrak{s} (\gamma_{2^m + 1}) \subseteq K^{(\times 2^m)}$ (as noted above) and:
\begin{center}
$[K^{(\times 2^{m-1})},\gamma_{2^m + 1}] \subseteq K^{(\times 2^m)}$
\end{center}
by Corollary \ref{degreecorollKs} (applied to both $K^{(\times 2^{m-1})}$ and $K^{(\times 2^m)}$). 
Thus $\overline{\mathfrak{s}}_m$ is indeed well-defined. 
\end{proof}

We now come to the heart of the proof of our diameter bound: 
using Propositions \ref{commspropnGrig} and \ref{squaring}, 
we show that if a symmetric subset $X \subseteq \mathfrak{G}$ 
contains an approximation to every element of $\mathfrak{G}$ 
up to an error lying in $\gamma_{2^m + 1}$, 
then every element of $\mathfrak{G}$ is approximated, up to 
an error in the (much smaller) subgroup $\gamma_{2^{m+1} + 1}$, 
by a short word in $X$. 

\begin{propn} \label{keypropnGrig}
Let $m \geq 2$ and let $X \subseteq \Gamma$ be a symmetric subset 
such that: 
\begin{equation} \label{Grighypeqn}
X \gamma_{2^m + 1} = \mathfrak{G}\text{.} 
\end{equation}
Then: 
\begin{equation} \label{Grigconcleqn}
X^{35} \gamma_{2^{m+1} + 1}  = \mathfrak{G}\text{.}
\end{equation}
\end{propn}

\begin{proof}
First, note that Proposition \ref{squaring} immediately implies: 
\begin{equation} \label{Grigsqringeqn}
X^2 K^{(\times 2^m)} \supseteq \gamma_{2^m + 1}\text{.}
\end{equation}

Second, we combine (\ref{Grighypeqn}) 
and Proposition \ref{commspropnGrig} (ii) 
with Corollary \ref{SKlem2} to conclude: 
\begin{equation} \label{Grigeqn2}
K^{(\times 2^m)} \subseteq X^4 \gamma_{2^m + 2^{m-1} + 1}\text{.}
\end{equation}
Taking stock of what we have thus far, (\ref{Grighypeqn}), (\ref{Grigsqringeqn}) and (\ref{Grigeqn2}) combine to give: 
\begin{equation} \label{Grigeqn3}
X^7 \gamma_{2^m + 2^{m-1} + 1} = \mathfrak{G} \text{.}
\end{equation}
Finally, we combine (\ref{Grigeqn3}) and (\ref{Grigeqn4}) 
with Corollary \ref{SKlem2} and conclude: 
\begin{equation} \label{Grigeqn5}
\gamma_{2^m + 2^{m-1} + 1} \subseteq X^{28} \gamma_{2^{m+1} + 1}\text{.}
\end{equation}
The required conclusion (\ref{Grigconcleqn}) is now immediate from (\ref{Grigeqn3}) and (\ref{Grigeqn5}). 
\end{proof}

\begin{proof}[Proof of Theorem \ref{mainthmGrig}]
Let $S \subseteq \mathfrak{G} / \gamma_n (\mathfrak{G})$ be a generating set. 
If $n \leq 5$ then: 
\begin{center}
$\diam (\mathfrak{G} / \gamma_n (\mathfrak{G}),S) 
\leq \lvert \mathfrak{G} : \gamma_n (\mathfrak{G}) \rvert$ 
\end{center}
is bounded by an absolute constant $\tilde{C}$. 
Otherwise let $\tilde{S} \subseteq \mathfrak{G}$ be any subset whose image in 
$\mathfrak{G} / \gamma_n (\mathfrak{G})$ is $S$, 
and let $m \in \mathbb{N}$ be such that $2^{m-1} + 1 < n \leq 2^m + 1$. 
Then $B_{\tilde{S}} (\tilde{C}) \gamma_5 (\mathfrak{G}) = \mathfrak{G}$, 
and by repeated application of Proposition \ref{keypropnGrig}, 
\begin{center}
$B_{\tilde{S}} (35^{m-2} \tilde{C}) \gamma_{2^m + 1} (\mathfrak{G}) = \mathfrak{G}$
\end{center}
so that $\diam (\mathfrak{G} / \gamma_n (\mathfrak{G}),S) \leq 35^{m-2} \tilde{C} \ll n^{\log (35)/\log (2)}$. 

The result now follows from Theorem \ref{centralindexGrig}. 
\end{proof}

\begin{rmrk}
Note that the above proof facilitates straightforward 
computation of the implied constant from the statement of 
Theorem \ref{mainthmGrig}. 
\end{rmrk}

\begin{proof}[Proof of Theorem \ref{stabdiamcorollGrig}]
By Lemma \ref{Kinstablem} and Corollary \ref{degreecorollKs}, 
\begin{center}
$\gamma_{2^{n+1}+1} (\mathfrak{G}) \leq \Stab_{\mathfrak{G}} (n+1)$
\end{center}
so: 
\begin{center}
$\diam (\mathfrak{G}/\Stab_{\mathfrak{G}} (n+1)) 
\leq \diam (\mathfrak{G}/\gamma_{2^{n+1}+1} (\mathfrak{G}))$. 
\end{center}
The result is now immediate from Theorem \ref{mainthmGrig}, 
Lemma \ref{stabindexlem} and Theorem \ref{centralindexGrig}. 
\end{proof}

\begin{rmrk} \normalfont
It is very likely that detailed knowledge of the lower central 
series of $\mathfrak{G}$ is not required to prove 
Theorem \ref{stabdiamcorollGrig}. 
Rather, one could give a direct proof of a diameter bound 
for $\mathfrak{G} / K^{(\times 2^n)}$ 
and deduce Theorem \ref{stabdiamcorollGrig} from this, 
much as we do for the Gupta-Sidki groups in the following Section. 
The reason for organizing the proof of Theorem 
\ref{stabdiamcorollGrig} as it appears here is historical. 
Theorem \ref{mainthmGrig} was the first of our results to be proved,  
followed be a direct proof of Theorem \ref{centralcorollGS}, 
using the results of \cite{Bar} 
on the lower central series of $\Gamma_{(3)}$. 
The case $p=3$ of Theorem \ref{stabcorollGS} was then deduced 
from Theorem \ref{centralcorollGS} 
(much as Theorem \ref{stabdiamcorollGrig} follows from 
Theorem \ref{mainthmGrig}). 
The (arguably more natural) proof of Theorem \ref{stabcorollGS} 
from Theorem \ref{mainthmGS} was a later response 
to the need to avoid assuming knowledge of the lower central 
series of $\Gamma_{(p)}$ in proving Theorem \ref{stabcorollGS} 
for higher $p$ 
(to the author's knowledge, the lower central series of 
$\Gamma_{(p)}$ has not been computed for $p \geq 5$). 
\end{rmrk}

\begin{rmrk} \normalfont
Proposition \ref{squaring} hints at the possibility 
of an alternative approach to proving diameter bounds 
for sequences of groups, 
following the same broad lines as the Solovay-Kitaev procedure 
employed here and in \cite{Brad}, but using power-words 
instead of commutator words. 
This will be explored further elsewhere \cite{Brad1}. 
\end{rmrk}

\section{Proofs for the Gupta-Sidki Groups} \label{GSproofsection}

In this section we prove Theorem \ref{mainthmGS} 
and, from it, Theorem \ref{stabcorollGS}. 
For the remainder of the section we shall write 
$\Gamma$ for $\Gamma_{(p)}$. 
Any assumptions on the prime $p$ in what follows 
will be made explicit in the appropriate place. 

The following notation will be useful in the sequel: 
for any $g \in \Aut (\mathcal{T}_{\mathcal{A}})$, \linebreak
let $\mathbf{0}(g) = (g,1,\ldots,1) \in \Aut (\mathcal{T}_{\mathcal{A}})^{(\times p)}$, 
and for $0 \leq j \leq p-1$, let \linebreak$(\mathbf{j+1})(g) = [a,\mathbf{j}(g)]$. 
Hence for $0 \leq j \leq p-1$, 
\begin{center}
$\mathbf{j}(g)_i = \big\lbrace \begin{array}{cc} g^{\alpha_{j,i}} & 1 \leq i \leq j+1 \\ 1 & \text{otherwise} \end{array}$, where $\alpha_{j,i} = (-1)^{i+1} \left( \begin{array}{c} j \\ i-1 \end{array} \right)$. 
\end{center}
Note that $\alpha_{p-1,i} \equiv 1 \mod p$ for all $1 \leq i \leq p$. 
It follows that: 
\begin{center}
$(\mathbf{p-1})(a) = (a,\ldots,a)$, $(\mathbf{p-1})(b) = (b,\ldots,b)$ 
\end{center}
(since $a,b$ have order $p$) and:
\begin{center}
$(\mathbf{p-1})(g) \equiv (g,\ldots,g) \mod K^{(\times p^2)}$ 
for any $g \in K$ 
\end{center}
(since $K \leq \Stab_{\Gamma} (1)$ and $\Gamma/K$ has exponent $p$, 
so too does $K / K^{(\times p)}$). 

Note also that in this notation, 
\begin{center}
$b = \mathbf{1}(a) \mathbf{0}(b)^{a^{-1}} 
= \mathbf{0}(b)^{a^{-1}} \mathbf{1}(a)$. 
\end{center}
It follows from the definition of the elements $x_i$ 
that for $1 \leq i \leq p-1$: 
\begin{equation} \label{constvecapprox}
x_i \big( (\mathbf{i+1})(a) \mathbf{i}(b)^{a^{-1}} \big)^{-1} \in K^{(\times p)}
\end{equation}
and in particular:
\begin{equation} \label{constvecapproxa}
x_{p-2} (a,\ldots,a)^{-1} \in (\Stab_{\Gamma} (1)) ^{(\times p)}
\end{equation}

\begin{equation} \label{constvecapproxb}
x_{p-1} (b,\ldots,b)^{-1} \in K^{(\times p)}
\end{equation}
These observations facilitate the completion of:

\begin{proof}[Proof of Corollary \ref{GSbranchprop} (ii)]
Since $K \leq \Stab_{\Gamma} (1)$, 
$K / K^{(\times p)}$ is an elementary abelian $p$-group. 
Moreover, from Proposition \ref{Garridoprop} (i) and (iii) we have 
\linebreak$\lvert K : K^{(\times p)} \rvert = p^{p-1}$, 
so it suffices to check that the images of $x_1 , \ldots , x_{p-1}$ 
in $K / K^{(\times p)}$ are linearly independent. 

Embedding $\Stab_{\Gamma} (1)$ into $\Gamma^{(\times p)}$ 
(via $\psi$), we have an induced embedding of $K / K^{(\times p)}$ 
into $(\Gamma/K)^{(\times p)}$. 
By Proposition \ref{Garridoprop} (i), 
\begin{center}
$(\Gamma/K)^{(\times p)} = \langle Ka , Kb \rangle^{(\times p)}
\cong C_p ^{(\times 2p)}$. 
\end{center}
From (\ref{constvecapprox}), the image of $x_i$ in 
$(\Gamma/K)^{(\times p)}$ is 
$v_i = (\mathbf{i+1})(Ka) \mathbf{i}(Kb)^{a^{-1}}$. 
For all $1 \leq i \leq p-2$, $v_i$ has non-zero $Ka$-component 
in its $(i+2)$th entry, but zero $Ka$-component in its $j$th entry for all $j \geq i+3$, so that $v_i$ is linearly independent of 
$\langle v_1 , \ldots , v_{i-1} \rangle$. 

Similarly, for all $1 \leq i \leq p-2$, 
$v_i$ has zero $Kb$-component in its $(p-1)$th entry, 
whereas $v_{p-1}$ has non-zero $Kb$-component there, 
so $v_{p-1}$ is independent of $\langle v_1, \ldots ,v_{p-2} \rangle$. 
\end{proof}

\begin{coroll} \label{GSStabindexcoroll}
For all $n \geq 1$, 
\begin{center}
$\lvert \Gamma : \Stab_{\Gamma} (n) \rvert 
\geq p^{(p-2)(p^{n-1}-1)+1}$. 
\end{center}
\end{coroll}

\begin{proof}
It suffices to check 
$\lvert \Stab_{\Gamma} (m) : \Stab_{\Gamma} (m+1) \rvert
\geq p^{(p-2)(p^{m-1}-1)}$ for all $m \geq 1$. 

Embedding $\Stab_{\Gamma} (m)$ into $\Gamma^{(\times p^m)}$ 
(via repeated application of $\psi$) we have an induced embedding:
\begin{center}
$\Stab_{\Gamma} (m) / \Stab_{\Gamma} (m+1) 
\hookrightarrow (\Gamma / \Stab_{\Gamma} (1))^{(\times p^m)}
\cong C_p ^{(\times p^m)} $
\end{center}
(with each $(\Gamma / \Stab_{\Gamma} (1))$-factor generated by 
$\Stab_{\Gamma} (1)a$). 
From (\ref{constvecapprox}), the image of 
$x_i \in \Stab_{\Gamma} (1)$ in 
$(\Gamma/\Stab_{\Gamma} (1))^{(\times p)}$ is 
$(\mathbf{i+1})(\Stab_{\Gamma} (1)a)$. 

Arguing as in the proof of Corollary \ref{GSbranchprop} (ii), 
the elements: 
\begin{equation*}
\Big( \prod_{j=1} ^{p-2} x_j ^{\lambda_{i,j}} \Big)_{i=1} ^{p^m}
\end{equation*}
are distinct modulo $\Stab_{\Gamma} (m+1)$ 
as the $p^{m-1} \times (p-2)$ coefficients $\lambda_{i,j}$ vary over 
$\lbrace 0,1,\ldots,p-1\rbrace$, 
and the required result follows. 
\end{proof}

We also introduce some further normal subgroups which will be useful 
``placeholders'' for our induction in the proof of Theorem \ref{mainthmGS}. 
For $1 \leq i \leq p$, let: 
\begin{center}
$L_i = \big\langle x_i , \ldots , x_{p-1} , K^{(\times p)} \big\rangle \leq K$
\end{center}
(with $L_p = K^{(\times p)}$ by convention). 
For $0 \leq i \leq p$, let: 
\begin{equation*}
\begin{aligned}
K_{\mathbf{i}} ^{(\times p)}
 &=  \big\langle \mathbf{i}(x_1) , \ldots , (\mathbf{p-1})(x_1), L_2 ^{(\times p)} \big\rangle \leq K^{(\times p)}\text{;} \\K_{\mathbf{i}} ^{(\times p^{n+1})} &= (K_{\mathbf{i}} ^{(\times p)}) ^{(\times p^n)} 
 \end{aligned}
\end{equation*}
(with the convention $K_{\mathbf{p}} ^{(\times p)} = L_2 ^{(\times p)}$). 
We thus have descending chains of subgroups:
\begin{center}
$L_2 ^{(\times p^n)} = K_{\mathbf{p}} ^{(\times p^n)} \leq K_{\mathbf{p-1}} ^{(\times p^n)} 
\leq \ldots \leq K_{\mathbf{1}} ^{(\times p^n)} \leq 
K_{\mathbf{0}} ^{(\times p^n)} = K^{(\times p^n)}$ 
\end{center}
and:
\begin{center}
$K^{(\times p^{n+1})} = L_{p} ^{(\times p^n)} \leq L_{p-1} ^{(\times p^n)} \leq \ldots \leq L_2 ^{(\times p^n)}$. 
\end{center}

\begin{lem}
The following are normal in $\Gamma_{(p)}$, for all $n \in \mathbb{N}$: 
\begin{itemize}
\item[(i)] $K^{(\times p^n)}$; 

\item[(ii)] $L_i ^{(\times p^n)}$, for $1 \leq i \leq p-1$; 

\item[(iii)] $K_{\mathbf{i}} ^{(\times p^{n+1})}$, for $0 \leq i \leq p-1$. 

\end{itemize}
\end{lem}

\begin{proof}
First, by induction we reduce to the case $n=0$. 
For let $H = K^{(\times p^n)}$, \linebreak$K_{\mathbf{i}} ^{(\times p^{n+1})}$ or $L_i ^{(\times p^n)}$. 
If $H \vartriangleleft \Gamma_{(p)}$, then $H^{(\times p)}$ is normalised by $a$ 
(which permutes the factors) and by $b$ (which acts on each factor as $a^{\pm 1}$ or $b$). 

Now certainly $K \vartriangleleft \Gamma_{(p)}$, so we have (i). 
In the other cases, 
we check that conjugates of a generating set for the subgroup, by the generators of $\Gamma_{(p)}$, 
lie in the subgroup. 

For (ii), consider the conjugation action of $\Gamma_{(p)}$ on $K / K^{(\times p)}$. 
We have \linebreak$x_i ^a = x_i x_{i+1} ^{-1}$ for $1 \leq i \leq p-2$; 
$x_{p-1} ^a \equiv x_{p-1} \mod K^{(\times p)}$ 
and $x_i ^b \equiv x_i  \mod K^{(\times p)}$ for $1 \leq i \leq p-1$. 
Thus $L_i \vartriangleleft \Gamma_{(p)}$. 

For (iii), consider the conjugation action of $\Gamma_{(p)}$ on $K^{(\times p)} / L_2 ^{(\times p)}$. \linebreak
We have $\mathbf{i}(x_1)^a = \mathbf{i}(x_1) (\mathbf{i+1})(x_1)^{-1}$ for $1 \leq i \leq p-2$; \linebreak
$(\mathbf{p-1})(x_1)^a \equiv (\mathbf{p-1})(x_1) \mod L_2 ^{(\times p)}$ 
and $\mathbf{i}(x_1)^b \equiv \mathbf{i}(x_1)$ for $1 \leq i \leq p-1$. 
Thus $K_{\mathbf{i}} ^{(\times p)} \vartriangleleft \Gamma_{(p)}$. 
\end{proof}

To apply Lemma \ref{SKlem1} in our induction, 
we will need certain commutators of elements in $\Gamma$ 
to lie in a sufficiently deep subgroup. 
As such, we note the following. 

\begin{lem} \label{subgrpincllem}
\begin{itemize}
\item[(i)] $[\Stab_{\Gamma} (1),\Stab_{\Gamma} (1)] \leq K^{(\times p)}$; 
\item[(ii)] $[\Stab_{\Gamma} (1),K^{(\times p)}] \leq L_2 ^{(\times p)}$. 
\end{itemize}
\end{lem}

\begin{proof}
For (i), note that any element of $\Stab_{\Gamma} (1)$ is a $p$-tuple of elements of $\Gamma$, 
so any commutator of such is a $p$-tuple of elements of $[\Gamma,\Gamma] = K$. 

Likewise in (ii), every element of $[\Stab_{\Gamma} (1),K^{(\times p)}]$ 
is a $p$-tuple of elements of $[\Gamma,K]$. $K / L_2$ is generated by $x_1$, 
and $\Gamma$ is generated by $a$ and $b$, so, since $L_2 \vartriangleleft K$, 
it suffices to check that $[a,x_1] , [b,x_1] \in L_2$ (the former by definition; the latter by (i)). 
\end{proof}

We now describe our approximations to elements of $\Gamma$ by commutators. 
These are split between the next four propositions, 
which are closely reminiscent of Proposition \ref{commspropnGrig} 
in our proof for Grigorchuk's group. 
The first of these allows us to approximate elements of 
$K ^{(\times p^2)}$ up to an error in $K_{\mathbf{1}} ^{(\times p^2)}$. We remark that this is the only point in the proof of Theorem \ref{mainthmGS} at which different arguments are required according to the value of $p$. To be precise, we exhibit one construction of an approximation by commutators which is valid for $p \geq 7$; 
one for $p=5$ and one for $p=3$. 

Let $\mathfrak{c} : \Gamma \times \Gamma \rightarrow \Gamma$ be given by $\mathfrak{c}(g,h)=[g,h]$. 

\begin{propn} \label{commspropnGS2}
The restriction of $\mathfrak{c}$ to $L_{p-1} \times K ^{(\times p)}$ 
descends to a well-defined map:
\begin{center}
$\bar{\mathfrak{c}}_0 : (L_{p-1} / K^{(\times p^2)}) 
\times (K ^{(\times p)} / K ^{(\times p^2)}) 
\rightarrow \Gamma / K_{\mathbf{1}} ^{(\times p^2)}$ 
\end{center}
whose image contains $K ^{(\times p^2)} / K_{\mathbf{1}} ^{(\times p^2)}$. 
\end{propn}

\begin{proof}
Observe that by Lemma \ref{SKlem1}, 
for $g_1 \in L_{p-1}$, $h_1 \in K ^{(\times p)}$
and $g_2$, \linebreak$h_2 \in K^{(\times p^2)}$
we have: 
\begin{center}
$[g_1 g_2 , h_1 h_2] [g_1,h_1]^{-1} 
\in \big[ L_{p-1} , K^{(\times p^2)} \big]$. 
\end{center}
But $L_{p-1}$ is generated by $K^{(\times p)}$ and $x_{p-1}$, 
both of which lie in $\Stab_{\Gamma} (1) ^{(\times p)}$ 
(the latter since $x_{p-1}$ is congruent 
modulo $K^{(\times p)}$ to $(b,\ldots ,b)$). 
Hence: 
\begin{center}
$\big[ L_{p-1} , K^{(\times p^2)} \big] 
\subseteq \big[ \Stab_{\Gamma}(1),K^{(\times p)} \big]^{(\times p)}
\subseteq L_2 ^{(p^2)}
\subseteq K_{\mathbf{1}} ^{(\times p^2)}$ 
(by Lemma \ref{subgrpincllem} (ii)). 
\end{center}
Thus $\bar{\mathfrak{c}}_0$ is indeed well-defined. 

We now establish that the image of $\bar{\mathfrak{c}}_0$ 
contains $K ^{(\times p^2)} / K_{\mathbf{1}} ^{(\times p^2)}$. 
First suppose $p \geq 7$. We have: 
\begin{equation*}
\begin{aligned}
x_1 ^a & = (b,b^{-1}a,a^{-2},a,1,\ldots,1)\\
x_1 ^{a^{-2}} & = (a,1,\ldots,1,b,b^{-1}a,a^{-2})
\end{aligned}
\end{equation*}
so that $\mathbf{0}(x_1) = [x_1 ^{a^{-2}},x_1 ^a]$.
Moreover for $\lambda \in \mathbb{N}$,
\begin{center}
$\mathbf{0}(x_1) ^{\lambda} 
\equiv [x_1 ^{a^{-2}},(x_1 ^a)^{\lambda}] 
\mod \big[ \big[ K,K \big] K \big] \leq L_2 ^{(\times p)} \leq K_{\mathbf{1}} ^{(\times p)}$
\end{center}
(by Lemma \ref{subgrpincllem}). 
 
Now, any element of $K ^{(\times p^2)} / K_{\mathbf{1}} ^{(\times p^2)}$ is represented by a vector 
$\big( \mathbf{0}(x_1)^{\lambda_j} \big)_{j=1} ^p$ 
for some $(\lambda_j)_{j=1} ^p \in \mathbb{F}_p ^p$.  By the above, 
\begin{center}
$\big( \mathbf{0}(x_1)^{\lambda_j} \big)_{j=1} ^p 
\equiv \big[ (x_1 ^{a^{-2}})_{j=1} ^p , ((x_1 ^a)^{\lambda_j})_{j=1} ^p \big]
\mod K_{\mathbf{1}} ^{(\times p^2)}$ 
\end{center}
and we are done. 

Now suppose $p = 5$. We have $x_1 ^a = (b,b^{-1}a,a^{-2},a,1)$. 
Hence: 
\begin{center}
$[x_1,x_1 ^a] = (x_1,aba^{-2}b^{-1}a,1,1,1)$. 
\end{center}
An easy direct computation yields 
$aba^{-2}b^{-1}a \equiv x_1 \mod L_2$. So:
\begin{center}
$[x_1,x_1 ^a] 
\equiv \mathbf{0} (x_1) ^2 \mod K_{\mathbf{1}} ^{(\times 5)}$. 
\end{center}
As before, any element of $K ^{(\times 5^2)} / K_{\mathbf{1}} ^{(\times 5^2)}$ is represented by: 
\begin{center}
$\big( \mathbf{0}(x_1)^{\lambda_j} \big)_{j=1} ^5 
\equiv \big[ (x_1)_{j=1} ^5 , ((x_1 ^a)^{3 \lambda_j})_{j=1} ^5 \big] \mod K_{\mathbf{1}} ^{(\times 5^2)}$ 
\end{center}
for some $(\lambda_j)_{j=1} ^5 \in \mathbb{F}_5 ^5$, as required. 

Finally suppose $p = 3$. Recall that:
\begin{center}
$b = (a,a^{-1},b)$ and 
$x_1 = [a,b] = (b^{-1}a,a,ab)$. 
\end{center}
Thus: 
\begin{center}
$[b,x_1] = ([a,b^{-1} a],1,[b,ab])$
\end{center}
We compute: 
\begin{align*}
[a,b^{-1} a] &= ([b,a]^{b^{-1}})^{a^{-1}} \\
& \equiv (x^{-1})^{a^{-1}} \mod K^{(\times 3)}\text{ (by Lemma \ref{subgrpincllem} (i))} \\
& \equiv x_1 ^{-1} \mod L_2
\end{align*}
and: 
\begin{align*}
[b,ab]& = [b,a][[b,a],b]\\
& \equiv x_1 ^{-1} \mod K^{(\times 3)}\text{ (by Lemma \ref{subgrpincllem} (i))}
\end{align*}
so that:
\begin{equation} \label{bx1eqn}
[b,x_1] \equiv (x_1 ^{-1},1,x_1 ^{-1}) \mod L_2 ^{(\times 3)}
\equiv \mathbf{0} (x_1) \mod K_{\mathbf{1}}^{(\times 3)}\text{.} 
\end{equation}
Now: 
\begin{align*}
x_2 = [a,x_1] &= (b^{-1} a^{-1} b^{-1} a , a^{-1} b a,b) \\
&= (bx_1 [x_1,b^{-1}],bx_1 ^{-1},b) \\
& \equiv (bx_1,bx_1 ^{-1},b) \mod K^{(\times 3^2)} 
\text{ (by Lemma \ref{subgrpincllem} (i))}
\end{align*}
so for $\lambda_1,\lambda_2,\lambda_3 \in \mathbb{N}$, 
\begin{center}
$[x_2,(x_1 ^{\lambda_1},x_1 ^{\lambda_2},x_1 ^{\lambda_3})]
\equiv \big( [bx_1 ,x_1 ^{\lambda_1}],[bx_1 ^{-1},x_1 ^{\lambda_2}],[b,x_1 ^{\lambda_3}] \big)
\mod K_{\mathbf{1}} ^{(\times 3^2)}$ \\

\end{center}
(by the well-definedness of $\overline{\mathfrak{c}}_0$) whereas for $\mu,\lambda \in \mathbb{N}$, 
\begin{align*}
[bx_1 ^{\mu},x_1 ^{\lambda}] &= [b,x_1 ^{\lambda}] [[b,x_1 ^{\lambda}],x_1 ^{\mu}] \\
& \equiv [b,x_1 ^{\lambda}] \mod L_2 ^{(\times 3)}\\
& \equiv [b,x_1]^{\lambda} \mod L_2 ^{(\times 3)}
\text{ (by Lemma \ref{subgrpincllem} (ii))}
\end{align*}
so by (\ref{bx1eqn}) we have: 
\begin{center}
$[x_2,(x_1 ^{\lambda_1},x_1 ^{\lambda_2},x_1 ^{\lambda_3})]
\equiv (\mathbf{0}(x_1)^{\lambda_1},\mathbf{0}(x_1)^{\lambda_2},\mathbf{0}(x_1)^{\lambda_3}) 
\mod L_2 ^{(\times 3^2)}$
\end{center}
and every element of 
$K ^{(\times p^2)} / K_{\mathbf{1}} ^{(\times p^2)}$
is represented by 
$(x_1 ^{\lambda_1},x_1 ^{\lambda_2},x_1 ^{\lambda_3})$ 
for some $\lambda_i \in \mathbb{F}_3$, as required. 
\end{proof}

Second, we construct an approximation to elements of $K_{\mathbf{i+1}} ^{(\times p^2)}$ up to an error lying in $K_{\mathbf{i+2}} ^{(\times p^2)}$. 

\begin{propn} \label{commspropnGS}
Let $0 \leq i \leq p-2$. 
The restriction of $\mathfrak{c}$ to $L_{p-2} \times K_{\mathbf{i}} ^{(\times p^2)}$ 
descends to a well-defined map:
\begin{center}
$\bar{\mathfrak{c}}_i : (L_{p-2} / K^{(\times p^2)}) 
\times (K_{\mathbf{i}} ^{(\times p^2)} / K_{\mathbf{i+1}} ^{(\times p^2)}) 
\rightarrow \Gamma / K_{\mathbf{i+2}} ^{(\times p^2)}$ 
\end{center}
whose image contains $K_{\mathbf{i+1}} ^{(\times p^2)} / K_{\mathbf{i+2}} ^{(\times p^2)}$. 
\end{propn}

\begin{proof}
We check first that $\bar{\mathfrak{c}}_i$ is well-defined. 
Observe that by Lemma \ref{SKlem1}, 
for $g_1 \in L_{p-2}$, 
$h_1 \in K_{\mathbf{i}} ^{(\times p^2)}$, 
$g_2 \in K^{(\times p^2)}$ and 
$h_2 \in K_{\mathbf{i+1}} ^{(\times p^2)}$, 
we have: 
\begin{center}
$[g_1 g_2 , h_1 h_2] [g_1,h_1]^{-1} \in \big[ L_{p-2},K_{\mathbf{i+1}} ^{(\times p^2)} \big] \big[ K_{\mathbf{i}} ^{(\times p^2)},K^{(\times p^2)} \big]$. 
\end{center}
It therefore suffices to check that: 
\begin{equation} \label{GSincl1}
\big[ L_{p-2},K_{\mathbf{i+1}} ^{(\times p^2)} \big] \big[ K_{\mathbf{i}} ^{(\times p^2)},K^{(\times p^2)} \big] 
\subseteq K_{\mathbf{i+2}} ^{(\times p^2)}\text{.}
\end{equation}
But $L_{p-2} = \langle x_{p-2} , L_{p-1} \rangle$ and 
$L_{p-1} \subseteq \Stab_{\Gamma} (2)$, 
so that by Lemma \ref{subgrpincllem} (i), 
\begin{center}
$\big [L_{p-1},K_{\mathbf{i+1}} ^{(\times p^2)} \big] 
\subseteq K ^{(\times p^3)} 
\subseteq K_{\mathbf{i+2}} ^{(\times p^2)}$
\end{center}
while $x_{p-2} (a,\ldots,a)^{-1} \in \Stab_{\Gamma} (1)^{(\times p)}$, so: 
\begin{center}
$\big[ \langle x_{p-2} \rangle , K_{\mathbf{i+1}} ^{(\times p^2)} \big] 
\subseteq \big[ \langle a \rangle , K_{\mathbf{i+1}} ^{(\times p)} \big]^{(\times p)} 
\big[ \Stab_{\Gamma} (1) , K_{\mathbf{i+1}} ^{(\times p)} \big]^{(\times p)} 
\subseteq K_{\mathbf{i+2}} ^{(\times p^2)}$ 
\end{center}
(using Lemma \ref{subgrpincllem} (ii)). 
Thus $\big[ L_{p-2},K_{\mathbf{i+1}} ^{(\times p^2)} \big] 
\subseteq K_{\mathbf{i+2}} ^{(\times p^2)}$. Also, 
$K_{\mathbf{i}} ^{(\times p^2)}$,\linebreak$K^{(\times p^2)} \subseteq \Stab_{\Gamma} (3)$, so by Lemma \ref{subgrpincllem} (i),
\begin{center}
$\big[ K_{\mathbf{i}} ^{(\times p^2)},K^{(\times p^2)} \big] \subseteq K^{(\times p^3)} \subseteq K_{\mathbf{i+2}} ^{(\times p^2)}$. 
\end{center}
Thus (\ref{GSincl1}) is indeed satisfied, 
and $\bar{\mathfrak{c}}_i$ is well-defined. 
\\ \\
We now check that the image of $\bar{\mathfrak{c}}_i$ contains $K_{\mathbf{i+1}} ^{(\times p^2)} / K_{\mathbf{i+2}} ^{(\times p^2)}$. 
First note that for any $\lambda \in \mathbb{N}$, 
\begin{align*}
(\mathbf{i+1})(x_1)^{\lambda} &= (\mathbf{i+1})(x_1 ^{\lambda})\\
& = [a,\mathbf{i}(x_1 ^{\lambda})]\\
& = [a,\mathbf{i}(x_1) ^{\lambda}]\text{.} 
\end{align*}
Now, every element of $K_{\mathbf{i+1}} ^{(\times p^2)} / K_{\mathbf{i+2}} ^{(\times p^2)}$ is represented by an element:
\begin{equation} \label{exprniplus1x1}
\big( (\mathbf{i+1})(x_1)^{\lambda_j} \big)_{j=1} ^p
 = \big[ (a)_{j=1} ^p , (\mathbf{i}(x_1) ^{\lambda_j})_{j=1} ^p \big] 
\end{equation}
for some $(\lambda_j)_{j=1} ^p \in \mathbb{N} ^p$.
From (\ref{constvecapproxa}), there exist 
$y_1 , \ldots , y_p \in \Stab_{\Gamma} (1)$ such that:
\begin{center}
$x_{p-2} = (a y_j)_{j=1} ^p$. 
\end{center}
For $1 \leq j \leq p$,
\begin{center}
$[ay_j,\mathbf{i}(x_1) ^{\lambda_j}] 
[a,\mathbf{i}(x_1) ^{\lambda_j}]^{-1} 
\in \big[ \Stab_{\Gamma} (1) , K_{\mathbf{i}} ^{(\times p)} \big] 
\subseteq K_{\mathbf{i+2}} ^{(\times p)}$
\end{center}
by Lemma \ref{SKlem1} and Lemma \ref{subgrpincllem} (ii). 
Combining with (\ref{exprniplus1x1}) we have:
\begin{center}
$\big( (\mathbf{i+1})(x_1)^{\lambda_j} \big)_{j=1} ^p \equiv 
\mathfrak{c} \big( x_{p-2} , (\mathbf{i}(x_1) ^{\lambda_j})_{j=1} ^p \big) 
\mod K_{\mathbf{i+2}} ^{(\times p^2)}$. 
\end{center}
Since $x_{p-2} \in L_{p-2}$ and 
$(\mathbf{i}(x_1) ^{\lambda_j})_{j=1} ^p 
\in K_{\mathbf{i}} ^{(\times p^2)}$
we are done. 
\end{proof}

Finally, we construct an approximation to elements of $L_{i+1} ^{(\times p^2)}$ up to an error lying in $L_{i+2} ^{(\times p^2)}$. 

\begin{propn} \label{commspropnGS3}
Let $1 \leq i \leq p-2$. 
The restriction of $\mathfrak{c}$ to $L_{p-2} ^{(\times p)} \times L_i ^{(\times p^2)}$ 
descends to a well-defined map:
\begin{center}
$\bar{\bar{\mathfrak{c}}}_i:(L_{p-2} ^{(\times p)}/K ^{(\times p^2)}) 
\times (L_i ^{(\times p^2)} / L_{i+1} ^{(\times p^2)}) \rightarrow 
\Gamma / L_{i+2} ^{(\times p^2)}$ 
\end{center}
whose image contains $L_{i+1} ^{(\times p^2)} / L_{i+2} ^{(\times p^2)}$. 
\end{propn}

\begin{proof}
We check first that $\bar{\bar{\mathfrak{c}}}_i$ is well-defined. 
Observe that by Lemma \ref{SKlem1}, 
for $g_1 \in L_{p-2} ^{(\times p)}$, 
$h_1 \in L_i ^{(\times p^2)}$, 
$g_2 \in K ^{(\times p^2)}$ and 
$h_2 \in L_{i+1} ^{(\times p^2)}$, 
we have: 
\begin{center}
$[g_1 g_2 , h_1 h_2] [g_1,h_1]^{-1} \in \big[ L_{p-2} ^{(\times p)} , L_{i+1} ^{(\times p^2)} \big] 
\big[ K^{(\times p^2)} , L_i ^{(\times p^2)} \big]$. 
\end{center}
It therefore suffices to check that: 
\begin{equation} \label{GSincl2}
\big[ L_{p-2} ^{(\times p)} , L_{i+1} ^{(\times p^2)} \big] 
\big[ K^{(\times p^2)} , L_i ^{(\times p^2)} \big]
\subseteq L_{i+2} ^{(\times p^2)}\text{.}
\end{equation}
Certainly, $\big[ K^{(\times p^2)} , L_i ^{(\times p^2)} \big] 
\leq K^{(\times p^3)} \leq L_{i+2} ^{(\times p^2)}$ 
(by Lemma \ref{subgrpincllem} (i)). 
Meanwhile, 
\begin{center}
$\big[ L_{p-2} ^{(\times p)} , L_{i+1} ^{(\times p^2)} \big] 
\subseteq \big[ L_{p-2} , L_{i+1} ^{(\times p)} \big]^{(\times p)}$
\end{center}
so it suffices to check that 
$\big[ L_{p-2} , L_{i+1} ^{(\times p)} \big] 
\subseteq L_{i+2} ^{(\times p)}$. 

$L_{p-2}$ is generated by $x_{p-2}$ and 
$L_{p-1} \subseteq \Stab_{\Gamma} (1)^{(\times p)}$ 
(the latter inclusion holds because $L_{p-1}$ 
is generated by $K^{(\times p)}$ and 
$x_{p-1} \in (b,\ldots ,b) K^{(\times p)}$) 
Thus:
\begin{center}
$\big[ L_{p-1},L_{i+1} ^{(\times p)} \big] 
\subseteq \big[ \Stab_{\Gamma} (1),L_{i+1} \big]^{(\times p)} 
\leq K^{(\times p^2)} \leq L_{i+1} ^{(\times p)}$. 
\end{center}
Meanwhile, $L_{i+1} ^{(\times p)} / L_{i+2} ^{(\times p)}$ 
is generated by the images of 
$\mathbf{0} (x_{i+1}) , \ldots , \mathbf{(p-1)} (x_{i+1})$, 
so by (\ref{constvecapproxa}), 
\begin{center}
$\big[ \langle x_{p-2} \rangle , L_{i+1} ^{(\times p)} \big] 
\subseteq\big[\langle a\rangle,\langle x_{i+1}\rangle\big]^{(\times p)} L_{i+2} ^{(\times p)}
\subseteq L_{i+2} ^{(\times p)}$. 
\end{center}
Thus (\ref{GSincl2}) is indeed satisfied, 
and $\bar{\bar{\mathfrak{c}}}_i$ is well-defined. 
\\ \\
We now check that the image of $\bar{\bar{\mathfrak{c}}}_i$ contains $L_{i+1} ^{(\times p^2)} / L_{i+2} ^{(\times p^2)}$. 

First, for any $\lambda \in \mathbb{N}$, 
\begin{center}
$x_{i+1} ^{\lambda} = [a,x_i]^{\lambda} \equiv [a,x_i ^{\lambda}] \mod K^{(\times p)}$ (by Lemma \ref{subgrpincllem} (i)). 
\end{center}
As in the proof of Proposition \ref{commspropnGS}, 
there exist, by (\ref{constvecapproxa}), 
$y_1 , \ldots , y_p \in \Stab_{\Gamma} (1)$ such that:
\begin{center}
$x_{p-2} = (a y_j)_{j=1} ^p$. 
\end{center}
For $1 \leq j \leq p$,
\begin{center}
$[ay_j,x_i ^{\lambda_j}] 
[a,x_i ^{\lambda_j}]^{-1} 
\in  \big[ K , K \big] \leq K^{(\times p)}$
\end{center}
(by Lemma \ref{subgrpincllem} (i)). 
Thus for $\lambda_1 , \ldots , \lambda_p \in \mathbb{N}$, 
\begin{center}
$[x_{p-2},(x_i ^{\lambda_1},\ldots,x_i ^{\lambda_p})] \equiv (x_{i+1} ^{\lambda_1},\ldots,x_{i+1} ^{\lambda_p}) \mod K^{(\times p^2)}$. 
\end{center}
Now every element of $L_{i+1} ^{(\times p^2)} / L_{i+2} ^{(\times p^2)}$ is represented by: 
\begin{center}
$(x_{i+1}^{\lambda_1} , \ldots , x_{i+1}^{\lambda_{p^2}})$
\end{center}
for some $(\lambda_j)_{j=1} ^{p^2} \in \mathbb{N} ^{p^2}$. From the above,
\begin{center}
$(x_{i+1}^{\lambda_1} , \ldots , x_{i+1}^{\lambda_{p^2}}) 
\equiv \big[ (x_{p-2})_{j=1} ^p,(x_i ^{\lambda_j})_{j=1} ^{p^2} \big]
\mod K^{(\times p^3)} \leq L_{i+2} ^{(\times p^2)}$
\end{center}
and we have $(x_{p-2})_{j=1} ^p \in L_{p-2} ^{(\times p)}$, 
$(x_i ^{\lambda_j})_{j=1} ^{p^2} \in L_i ^{(\times p^2)}$ 
so the result follows.  
\end{proof}

Up to now, we have concentrated on approximating, by commutators, 
elements lying in $K^{(\times p^2)}$ but outside $K^{(\times p^3)}$. 
We can however quickly extend these approximations 
to elements lying between 
$K^{(\times p^m)}$ amd $K^{(\times p^{m+1})}$ 
for arbitrary $m \geq 2$. 
Indeed, since the conclusions of 
Propositions \ref{commspropnGS2}-\ref{commspropnGS3}
concern only computations within the group $K/K^{(\times p^3)}$
the generalisation from the case $m=2$ 
is immediate from the identification: 
\begin{center}
$K^{(\times p^{m-2})} / K^{(\times p^{m+1})} 
\cong (K / K^{(\times p^3)})^{(\times p^{m-2})}$
\end{center}
and the observation that, 
in a direct product of groups, 
a tuple of commutators of elements of the factors is the commutator 
of the tuples of those same elements. 

\begin{propn} \label{commspropnGS4}
Let $m \geq 2$. 
\begin{itemize}
\item[(i)] The restriction of $\mathfrak{c}$ to 
$L_{p-1} ^{(\times p^{m-2})} \times K ^{(\times p^{m-1})}$ 
descends to a well-defined map:
\begin{center}
$\bar{\mathfrak{c}}_{0,m} : (L_{p-1} ^{(\times p^{m-2})} / K^{(\times p^m)}) 
\times (K ^{(\times p^{m-1})} / K ^{(\times p^m)}) 
\rightarrow \Gamma / K_{\mathbf{1}} ^{(\times p^m)}$ 
\end{center}
whose image contains $K ^{(\times p^m)} / K_{\mathbf{1}} ^{(\times p^m)}$. 

\item[(ii)] Let $0 \leq i \leq p-2$. 
The restriction of $\mathfrak{c}$ to $L_{p-2} ^{(\times p^{m-2})} \times K_{\mathbf{i}} ^{(\times p^m)}$ 
descends to a well-defined map:
\begin{center}
$\bar{\mathfrak{c}}_{i,m} : (L_{p-2} ^{(\times p^{m-2})} / K^{(\times p^m)}) 
\times (K_{\mathbf{i}} ^{(\times p^m)} / K_{\mathbf{i+1}} ^{(\times p^m)}) 
\rightarrow \Gamma / K_{\mathbf{i+2}} ^{(\times p^m)}$ 
\end{center}
whose image contains $K_{\mathbf{i+1}} ^{(\times p^m)} / K_{\mathbf{i+2}} ^{(\times p^m)}$. 

\item[(iii)] Let $1 \leq i \leq p-2$. 
The restriction of $\mathfrak{c}$ to $L_{p-2} ^{(\times p^{m-1})} \times L_i ^{(\times p^{m-2})}$ 
descends to a well-defined map:
\begin{center}
$\bar{\bar{\mathfrak{c}}}_{i,m}:(L_{p-2} ^{(\times p^{m-1})}/K ^{(\times p^m)}) 
\times (L_i ^{(\times p^m)} / L_{i+1} ^{(\times p^m)}) \rightarrow 
\Gamma / L_{i+2} ^{(\times p^m)}$
\end{center}
whose image contains $L_{i+1} ^{(\times p^m)} / L_{i+2} ^{(\times p^m)}$. 
\end{itemize}
\end{propn}

\begin{proof}
By the preceding discussion, (i), (ii) and (iii) follow, 
respectively, from Propositions \ref{commspropnGS2}, 
\ref{commspropnGS} and \ref{commspropnGS3}. 
\end{proof}

We are now ready to put everything together, 
and use the approximations from Proposition \ref{commspropnGS4} 
(i)-(iii) to prove a result closely analogous to 
Proposition \ref{keypropnGrig}: namely, if a symmetric subset 
$X \subseteq \Gamma$ contains an approximation 
to every element of $\Gamma$ up to an error lying in 
$K^{(\times p^m)}$, then every element of $\Gamma$ is 
approximated up to an error lying in $K^{(\times p^{m+1})}$ 
by a short word in $X$. 

\begin{propn} \label{keypropnGS}
Let $C_p$ be as in Theorem \ref{stabcorollGS}. 
Let $m \geq 2$ and let $X \subseteq \Gamma$ be 
a symmetric subset such that: 
\begin{equation} \label{GShypeqn}
X K^{(\times p^m)} = \Gamma\text{.}
\end{equation}
Then: 
\begin{equation} \label{GSconcleqn}
X^{C_p} K^{(\times p^{m+1})} = \Gamma\text{.}
\end{equation}
\end{propn}

\begin{proof}
The first step shall be to show that: 
\begin{equation} \label{GSstep1goal}
K ^{(\times p^m)} \subseteq X^4 K_{\mathbf{1}} ^{(\times p^m)}\text{.}
\end{equation}
To this end let $k \in K ^{(\times p^m)}$. 
By Proposition \ref{commspropnGS4} (i), 
there exist $g \in L_{p-1} ^{(\times p^{m-2})}$, 
$h \in K ^{(\times p^{m-1})}$ such that: 
\begin{center}
$[g,h] \equiv k \mod K_{\mathbf{1}} ^{(\times p^m)}$. 
\end{center}
From (\ref{GShypeqn}), there exist $x_g , x_h \in X$ such that: 
\begin{center}
$x_g \equiv g , x_h \equiv h \mod K ^{(\times p^m)}$
\end{center}
so by the well-definedness of the map $\bar{\mathfrak{c}}_{0,m}$ 
from Proposition \ref{commspropnGS4} (i), 
\begin{center}
$X^4 \ni [x_g,x_h] \equiv k \mod K_{\mathbf{1}} ^{(\times p^m)}$ 
\end{center}
and we have (\ref{GSstep1goal}). 

Define the integer sequence $(a_n)_n$ recursively by $a_0 = 4$ 
and $a_n = 2 a_{n-1} + 2$ for $n \geq 1$. 
The second step of the proof shall be to show that: 
\begin{equation} \label{GSstep2goal}
K_{\mathbf{i}} ^{(\times p^m)} \subseteq X^{a_i} K_{\mathbf{i+1}} ^{(\times p^m)}
\end{equation}
for $0 \leq i \leq p-1$. This shall be achieved by induction on $i$, 
using Proposition \ref{commspropnGS4} (ii) at each stage 
(and the base case $i=0$ being provided by (\ref{GSstep1goal})). 

For let $1 \leq i \leq p-1$ and let 
$k \in K_{\mathbf{i}} ^{(\times p^m)}$. 
By Proposition \ref{commspropnGS4} (ii), 
there exist $g \in L_{p-2} ^{(\times p^{m-2})}$, 
$h \in K_{\mathbf{i-1}} ^{(\times p^{m})}$ such that: 
\begin{center}
$[g,h] \equiv k \mod K_{\mathbf{i+1}} ^{(\times p^m)}$. 
\end{center}
From (\ref{GShypeqn}) and the induction hypothesis, 
there exist $x_g \in X$, $x_h \in X^{a_{i-1}}$ such that: 
\begin{center}
$x_g \equiv g \mod K ^{(\times p^m)}$, 
$x_h \equiv h \mod K_{\mathbf{i}} ^{(\times p^m)}$
\end{center}
so by the well-definedness of the map $\bar{\mathfrak{c}}_{i,m}$ 
from Proposition \ref{commspropnGS4} (ii), 
\begin{center}
$X^{a_i} = X^{2a_{i-1} + 2} \ni [x_g,x_h] \equiv k \mod K_{\mathbf{i+1}} ^{(\times p^m)}$ 
\end{center}
and we have (\ref{GSstep2goal}). 

Define the integer sequence $(b_n)_n$ recursively by 
$b_1 = \sum_{n=0} ^{p-1} a_n$ and \linebreak$b_{n+1} = 2 b_n + 2$ 
for $n \geq 2$. 
Combining the inclusions (\ref{GSstep2goal}) 
for $i$ from $0$ to $p-1$, we have: 
\begin{equation} \label{GSphase1}
K^{(\times p^m)} = K_{\mathbf{0}} ^{(\times p^m)} 
\subseteq X^{b_1} K_{\mathbf{p}} ^{(\times p^m)} 
= X^{b_1} L_2 ^{(\times p^m)}\text{.} 
\end{equation}
Our third objective shall be to show that: 
\begin{equation} \label{GSstep3goal}
L_{i} ^{(\times p^m)} 
\subseteq X^{b_i} L_{i+1} ^{(\times p^m)}
\end{equation}
for $1 \leq i \leq p-1$. 
This again shall be by induction on $i$, 
using Proposition \ref{commspropnGS4} (iii), 
the base case $i=1$ being provided by (\ref{GSphase1}). 

Thus let $2 \leq i \leq p-1$ and let 
$k \in L_{i} ^{(\times p^m)}$. 
By Proposition \ref{commspropnGS4} (ii), 
there exist $g \in L_{p-2} ^{(\times p^{m-1})}$, 
$h \in L_{i-1} ^{(\times p^m)}$ such that: 
\begin{center}
$[g,h] \equiv k \mod L_{i+1} ^{(\times p^m)}$. 
\end{center}
By (\ref{GShypeqn}) and the induction hypothesis, 
there exist $x_g \in X$, $x_h \in X^{b_{i-1}}$ such that: 
\begin{center}
$x_g \equiv g \mod K ^{(\times p^m)}$, 
$x_h \equiv h \mod L_{i} ^{(\times p^m)}$
\end{center}
so by the well-definedness of the map $\bar{\bar{\mathfrak{c}}}_{i,m}$ 
from Proposition \ref{commspropnGS4} (iii), 
\begin{center}
$X^{b_i} = X^{2b_{i-1} + 2} \ni [x_g,x_h] 
\equiv k \mod L_{i+1} ^{(\times p^m)} $ 
\end{center}
as desired. 

Finally, set $C_p = 1 + \sum_{i=1} ^{p-1} b_i$. 
Expressing $a_n$ and $b_n$ in closed form, 
$C_p$ is as in the statement of Theorems \ref{stabcorollGS} 
and \ref{mainthmGS}. 
We combine the inclusions (\ref{GSstep3goal}) 
for $i$ from $1$ to $p-1$ to obtain: 
\begin{center}
$K^{(\times p^m)}  = L_1 ^{(\times p^m)} 
\subseteq X^{C_p - 1} L_p ^{(\times p^m)} 
= X^{C_p - 1} K^{(\times p^{m+1})}$. 
\end{center}
Combining this last inclusion with (\ref{GShypeqn}), 
we have $\Gamma = X^{C_p} K^{(\times p^{m+1})}$, 
as required. 
\end{proof}

\begin{proof}[Proof of Theorem \ref{mainthmGS}]
Let $S \subseteq \Gamma / K^{(\times p^n)}$. 
If $n \leq 2$ then: 
\begin{center}
$\diam (\Gamma / K^{(\times p^n)}) 
\leq \lvert G : K^{(\times p^2)} \rvert = p^{p+1}$. 
\end{center}
If $n \geq 3$ then let $\tilde{S} \subseteq \Gamma$ 
be any subset whose image in $\Gamma / K^{(\times p^n)}$ is $S$. 
Then $B_{\tilde{S}} (p^{p+1}) K^{(\times p^2)} = \Gamma$, 
and by repeated application of Proposition \ref{keypropnGS}, 
\begin{center}
$B_{\tilde{S}} (p^{p+1} C_p ^{n-2}) K^{(\times p^n)} = \Gamma$. 
\end{center}
Thus: 
\begin{center}
$\diam (\Gamma / K^{(\times p^n)} , S) 
\leq p^{p+1} C_p ^{n-2} \ll_p p^{\log (C_p) n / \log (p)}$
\end{center}
The result follows, since 
$\lvert \Gamma : K^{(\times p^n)} \rvert = p^{p^n + 1}$. 
\end{proof}

\begin{proof}[Proof of Theorem \ref{stabcorollGS}]
By Lemma \ref{GSKStabincl}, we have: 
\begin{center}
$\diam (\Gamma / \Stab_{\Gamma} (n)) \leq \diam (\Gamma / K^{(\times p^{n-1})})$. 
\end{center}
The result now follows from 
$\lvert \Gamma : K^{(\times p^n)} \rvert = p^{p^n + 1}$; 
Theorem \ref{mainthmGS} and Corollary \ref{GSStabindexcoroll}. 
\end{proof}

\section{The Gupta-Sidki $3$-Group} \label{GS3sect}

In this Section we deduce Theorem \ref{centralcorollGS} 
from Theorem \ref{mainthmGS}. We shall require some facts about the lower central series of $\Gamma = \Gamma_{(3)}$, which were established in \cite{Bar}. 
Define the integer sequence $(\alpha_n)$ by $\alpha_1 = 1$, $\alpha_2 = 2$, 
$\alpha_n = 2 \alpha_{n-1} + \alpha_{n-2}$ for $n \geq 3$, and set $\beta_n = \sum_{i=1} ^n \alpha_i$. We have: 
\begin{align*}
\alpha_n =& \frac{1}{ 2 \sqrt{2}} \big((1 + \sqrt{2})^n - (1 - \sqrt{2})^n\big), \\
\beta_n =& \frac{1}{4} \big((1 + \sqrt{2})^{n+1} + (1 - \sqrt{2})^{n+1}-2\big)\text{.} 
\end{align*}

\begin{thm}[\cite{Bar}]
Let $X_1 , \ldots , X_n \in \lbrace \mathbf{0} , \mathbf{1} , \mathbf{2} \rbrace$. Then:
\begin{align*}
\deg X_1 \cdots X_n (x_1) &= 1 + \sum_{i=1} ^n X_i \alpha_i + \alpha_{n+1} \\
\deg X_1 \cdots X_n (x_2) &= 1 + \sum_{i=1} ^n X_i \alpha_i + 2 \alpha_{n+1}\text{.} 
\end{align*}
\end{thm}

\begin{coroll} \label{GS3inclcoroll}
For all $m \in \mathbb{N}$, 
\begin{center}
$K^{(\times 3^m)} \leq \gamma_{\alpha_{m+1}+1} (\Gamma) 
\leq \gamma_{\beta_m+1} (\Gamma)$.
\end{center}
\end{coroll}

\begin{coroll} \label{GS3ordercoroll}
For all $m \in \mathbb{N}$, 
\begin{center}
$\lvert \Gamma : \gamma_{\beta_m + 1} (\Gamma) \rvert 
= 3^{(3^m + 1)/2}$. 
\end{center}
\end{coroll}

\begin{proof}[Proof of Theorem \ref{centralcorollGS}]
For $n=1$ there is nothing to prove. 
Otherwise, let $m \in \mathbb{N}$ be such that 
$\beta_m + 1 \leq n \leq \beta_{m+1}+1$. 
Then by Corollary \ref{GS3inclcoroll}, 
\begin{align*}
\diam (\Gamma / \gamma_n (\Gamma)) 
& \leq \diam (\Gamma / K^{(\times 3^{m+1})}) \\
& \ll 3^{\log(111) m / \log (3)}\text{ (by Theorem \ref{mainthmGS})}\\
& \ll n^{\log (111)/\log(1+\sqrt{2})} \\
& \ll (\log \lvert \Gamma : \gamma_{\beta_m + 1} (\Gamma) \rvert)^{\log(111) / \log (3)}\text{ (by Corollary \ref{GS3ordercoroll})} \\
& \leq (\log \lvert \Gamma : \gamma_n (\Gamma) \rvert)^{\log(111) / \log (3)}\text{.} 
\end{align*}
\end{proof}

\section{Spectral Gap and Mixing Time} \label{gapsect}

For $G$ a finite group and $S \subseteq G$ a symmetric subset, 
let $A_S$ be the (normalized) adjacency operator on the Cayley 
graph $\Cay (G,S)$. $A_S$ is a self-adjoint operator of norm $1$; 
let its spectrum be: 
\begin{center}
$1 = \lambda_1 \geq \lambda_2 \geq\ldots\geq \lambda_{\lvert G \rvert} \geq -1$
\end{center}
with the eigenvalue $\lambda_1 = 1$ corresponding to the constant 
functionals on $G$. More generally, the $1$-eigenspace of $A_S$ 
is spanned by the indicator functions of the connected components of 
$\Cay (G,S)$; in particular $1 > \lambda_2$ if and only if $S$ 
generates $G$, and in this case the quantity $1 - \lambda_2$ 
is known as the \emph{spectral gap} of the pair $(G,S)$. 

The existence of a large spectral gap for a family 
of Cayley graphs is a matter of great interest. 
If a family of finite graphs of bounded valence with vertex sets of unbounded 
size possess a spectral gap bounded below by an absolute positive constant, 
then the graphs form an \emph{expander family}. 
Expander graphs (and especially expander \emph{Cayley} graphs) have 
multifarious applications across pure mathematics and 
theoretical computer science \cite{Lub}. 

Now let $\Gamma$ be $\mathfrak{G}$ or $\Gamma_{(p)}$ and let 
$(\Gamma_i)_i$ be one of the descending sequences 
of finite-index normal subgroups 
from Theorems \ref{stabdiamcorollGrig}-\ref{centralcorollGS}. 
Cayley graphs of the quotient groups $\Gamma/\Gamma_i$ 
do not in general form expander families: 
for instance if $S$ is a finite symmetric generating set for $\Gamma$, 
and $S_i$ is the image of $S$ in $\Gamma/\Gamma_i$, 
then the spectral gap of $\Cay (\Gamma/\Gamma_i , S_i)$ 
tends to $0$ as $i \rightarrow \infty$ 
(this follows from the fact that $\Gamma$ is amenable \cite{BaKaNe}
and $\Gamma_i$ exhausts $\Gamma$). 
We do however have a weaker lower bound on the spectral gap 
of any connected Cayley graph of $\Gamma/\Gamma_i$, 
coming from our upper bounds on diameter and the following general inequality. 

\begin{propn}[\cite{DiaSaCo} Corollary 3.1] \label{spectralgappropn}
Let $G$ be a finite group and let $S$ be a symmetric generating set. 
Then the spectral gap of $(G,S)$ is $\geq (\lvert S \rvert \diam (G,S)^2)^{-1}$.  
\end{propn}

Theorems \ref{stabdiamcorollGrig}-\ref{centralcorollGS} combine 
with Proposition \ref{spectralgappropn} to yield the following 
bounds on spectral gaps. 

\begin{coroll} \label{branchsgcoroll}
Let $S$ be an arbitrary generating set for the finite group $G$. 
Denote by $\epsilon (G,S)$ the spectral gap of the pair $(G,S)$. 
Let $C_p$ be as in Theorem \ref{stabcorollGS}. 
\begin{itemize}
\item[(i)] If $G = \mathfrak{G}/\Stab_{\mathfrak{G}} (n)$ then $\epsilon(G,S) = \Omega \big( \lvert S \rvert^{-1} \exp (-2\log(35)n) \big)$; 

\item[(ii)] If $G = \mathfrak{G}/\gamma_n(\mathfrak{G})$ then $\epsilon(G,S) = \Omega \big( \lvert S \rvert^{-1} n^{-2\log(35)/\log(2)} \big)$; 

\item[(iii)] If $G = \Gamma_{(p)}/\Stab_{\Gamma_{(p)}} (n)$ then $\epsilon(G,S) = \Omega_p \big( \lvert S \rvert^{-1}  \exp (-2\log (C_p)n) \big)$; 

\item[(iv)] If $G = \Gamma_{(p)}/K^{(\times p^n)}$ then $\epsilon(G,S) = \Omega_p \big( \lvert S \rvert^{-1}   \exp (-2\log (C_p)n) \big)$; 

\item[(v)] If $G = \Gamma_{(3)}/\gamma_n (\Gamma_{(3)})$ then $\epsilon(G,S) = \Omega \big( \lvert S \rvert^{-1} n^{-2 \log(111)/\log(1+\sqrt{2})}  \big)$. 

\end{itemize}
\end{coroll}

A second closely related numerical invariant of finite Cayley graphs 
is the \emph{mixing time}. This is a measure of the time taken 
for a lazy random walk on the Cayley graph to approach the uniform distribution. 
It may be defined as follows. Let $f_0 = \delta_e$ be the Dirac mass at 
the identity of $G$. The lazy random walk on $\Cay(G,S)$ 
is defined by the operator $T_S = (A_S + I)/2$, where $I$ is the 
identity operator on $\Cay(G,S)$, 
and describes the progress on $\Cay(G,S)$ 
of a particle which starts at the identity, 
and which at each step with equal probability either traverses an edge 
(chosen uniformly at random) or remains stationary. 
Recursively define $f_{l+1} = T_S (f_l)$, 
the distribution of the walk at time $l$. 
We may consider the walk to be well-mixed when $f_l$ is close 
to the uniform distribution, in some appropriate norm 
on the complex functionals on $G$. 
Here we focus on mixing with respect to the $\ell^{\infty}$-norm. 

\begin{defn}
Let $G$ be a finite group and $S$ be a symmetric generating set. 
The \emph{$\ell^{\infty}$-mixing time} of the pair $(G,S)$ 
is the smallest positive integer $l$ such that: 
\begin{center}
$\big\lVert f_l - \frac{1}{\lvert G \rvert} \chi_G \big\rVert_{\infty} \leq \frac{1}{2 \lvert G \rvert}$. 
\end{center}
\end{defn}

It may be easily seen that the LHS of the above inequality is 
a non-increasing function of $l$, 
so that once the random walk reaches its mixing time, 
it remains well-mixed thereafter. 
There is a close relationship between mixing time and spectral gap. 

\begin{propn}[\cite{Lova} Theorem 5.1]
Suppose the pair $(G,S)$ has spectral gap $\epsilon > 0$. 
Then there exists an absolute constant $C>0$ such that 
the $\ell^{\infty}$-mixing time of $(G,S)$ is at most 
$(C/\epsilon) \log \lvert G \rvert$. 
\end{propn}

Applying Proposition  to the conclusions of Corollary \ref{branchsgcoroll}, 
we have corresponding bounds on mixing times, as follows. 

\begin{coroll}
Denote by $\mu (G,S)$ the $\ell^{\infty}$-mixing time of the pair $(G,S)$. 
\begin{itemize}
\item[(i)] If $G = \mathfrak{G}/\Stab_{\mathfrak{G}} (n)$ then $\mu(G,S) = O \big( \lvert S \rvert \exp (\log(2450)n) \big)$; 

\item[(ii)] If $G = \mathfrak{G}/\gamma_n(\mathfrak{G})$ then $\mu(G,S) = O \big( \lvert S \rvert n^{2\log(35)/\log(2) + 1} \big)$; 

\item[(iii)] If $G = \Gamma_{(p)}/\Stab_{\Gamma_{(p)}} (n)$ then $\mu(G,S) = O_p \big( \lvert S \rvert \exp (\log (pC_p ^2)n) \big)$; 

\item[(iv)] If $G = \Gamma_{(p)}/K^{(\times p^n)}$ then $\mu(G,S) = O_p \big( \lvert S \rvert \exp (\log (pC_p ^2)n) \big)$; 

\item[(v)] If $G = \Gamma_{(3)}/\gamma_n (\Gamma_{(3)})$ then $\mu(G,S) = O \big( \lvert S \rvert n^{\log(36963)/\log(1+\sqrt{2})}  \big)$. 

\end{itemize}
\end{coroll}

\section{Growth in Branch Groups} \label{growthsect}

Given a finitely generated group $G$ and a finite generating set $S \subseteq G$, let $f_{(G,S)} (n) = \lvert B_S (n) \rvert$ be the \emph{growth function}. Although for a given group $G$, 
the function $f_{(G,S)}$ may vary according to the generating set $S$, 
it only does so up to an appropriate notion of equivalence of functions. 
As such, we may speak without ambiguity about groups of 
\emph{polynomial growth}, \emph{exponential growth} and so on 
(see \cite{dlHa} Chapters VI-VII). 

One of the key sources of interest in branch groups is the fact that they include many examples of groups with exotic growth behaviour. 
In particular, $\mathfrak{G}$ has \emph{intermediate growth}, 
that is: growth faster than any polynomial function but slower than any exponential function. 

The following elementary fact exhibits a relationship between growth and diameter. 

\begin{lem}
Let $F$ be a finite group, and let $\phi : G \rightarrow F$ 
be an epimorphism. 
Then: 
\begin{equation} \label{growthdiamineq}
f_{(G,S)} (\diam (F,\phi (S))) \geq \lvert F \rvert\text{.}
\end{equation}
\end{lem}

This inequality suggests the following definition, which is made by analogy with that of the diameter of a finite group. Let $f_G (n)$ be the minimal 
value of $f_{(G,S)} (n)$, as $S$ ranges over all finite generating subsets of $G$. From (\ref{growthdiamineq}) we immediately obtain: 
\begin{equation} \label{ugrowthdiamineq}
f_G (\diam (F)) \geq \lvert F \rvert\text{.}
\end{equation}
The relationship between growth and diameter 
can be exploited to yield information about both. 
For instance, from Theorems \ref{mainthmGrig} 
and \ref{mainthmGS} we have the following. 

\begin{coroll} \label{lbgrowthcoroll}
There exist constants $\alpha (p) > 0$ such that: 
\begin{center} 
$f_{\mathfrak{G}} (n) \gg \exp (\alpha (2) n ^{\beta(\mathfrak{G})})$ and $f_{\Gamma_{(p)}} (n) \gg_p \exp (\alpha (p) n ^{\beta(\Gamma_{(p)})})$, 
\end{center}
where $\beta(\mathfrak{G}) = \log (2)/\log(35) \approx 0.195$ 
and $\beta(\Gamma_{(p)}) = \log (p) / \log (C_2 (p))$ 
(here $C_2 (p)$ is as in Theorem \ref{stabcorollGS}). 
\end{coroll}

The bounds in Corollary \ref{lbgrowthcoroll} are not the best known: 
a slight modification of an argument of Grigorchuk \cite{Grig2} 
shows that if $G$ is any finitely generated residually virtually 
nilpotent group, then either $G$ is virtually nilpotent or 
$f_G (n) \geq \exp (\sqrt{n})$ (see \cite{Bar2}). 
In particular the latter conclusion applies to $\mathfrak{G}$ 
and $\Gamma_{(p)}$. 
It is however possible that improvements upon the diameter bounds in 
Theorems \ref{mainthmGrig} and \ref{mainthmGS} and their corollaries 
could yield new lower bounds on $f_{\mathfrak{G}}$ 
and $f_{\Gamma_{(p)}}$. 

Conversely, known upper bounds on the growth translate into 
lower bounds on the diameters of finite quotients. 
In the case of $\mathfrak{G}$, the best upper bound on the growth 
is the following result of Bartholdi. 

\begin{thm}[\cite{Grig1}] \label{ubgrowthGrig}
Let $a,b,c,d \in \mathfrak{G}$ be as in subsection \ref{Grigdefnsubsubsect} and let $S = \lbrace a,b,c,d \rbrace$. 
Then: 
\begin{center}
$f_{(\mathfrak{G},S)} (n) \ll \exp (n^\beta)$, 
\end{center}
where $\beta = \frac{\log (2)}{\log (2) - \log (\eta)} \approx 0.768$, 
for $\eta$ the real root of $X^3 + X^2 + X = 2$. 
\end{thm}

\begin{coroll} \label{lbdiamGrig}
There exist absolute constants $C , C^{\prime} > 0$ such that: 
\begin{center}
$\diam (\mathfrak{G} / \Stab_{\mathfrak{G}} (n)) 
\geq C(\log \lvert \mathfrak{G} : \Stab_{\mathfrak{G}} (n) \rvert)^{1/\beta}$; \\
$\diam (\mathfrak{G} / \gamma_n (\mathfrak{G})) 
\geq C^{\prime} (\log \lvert \mathfrak{G} : \gamma_n (\mathfrak{G}) \rvert )^{1/\beta}$
\end{center}
where $\beta$ is as in Theorem \ref{ubgrowthGrig}. 
\end{coroll}

Corollary \ref{lbdiamGrig} places a limit 
on the extent to which the constant $\log(35)/\log(2)$ 
appearing in Theorems \ref{stabdiamcorollGrig} 
and \ref{mainthmGrig} might be reduced 
(though it is almost certainly not sharp). 
It is unclear at this time whether the constant 
$1/\beta \approx 1.303$ in Corollary \ref{lbdiamGrig}
could be close to sharp. 

\subsection*{Acknowledgements}

I would like to thank Laurent Bartholdi and Alejandra Garrido for illuminating discussions. 
Parts of this work were completed while the author was a fellow-commoner at Trinity Hall, Cambridge. 
I would like to thank the fellows of the College for providing me with a warm welcome and pleasant working conditions.

\Addresses

\end{document}